\documentclass[journal]{IEEEtran}
\usepackage[full]{textcomp}
\usepackage{graphicx}
\usepackage{rotate}
\usepackage{times}
\usepackage{amsmath}
\usepackage{amscd}
\usepackage{amsthm}
\usepackage{amssymb}
\usepackage{epsfig}
\usepackage{color}
\usepackage{colortbl}
\usepackage[table]{xcolor}
\usepackage{url}
\usepackage{soul}
\usepackage{enumitem}
\usepackage{multirow}
\usepackage{calligra}
\usepackage{aurical}
\usepackage{xfrac}
\usepackage[T1]{fontenc}
\usepackage{ragged2e}
 \usepackage{mathtools,bbm}
\usepackage[cal=cm, frak=euler, scr=boondoxo]{mathalfa}
 \usepackage[rmdefault=cmr]{isomath}
 \usepackage{mathabx}
  \usepackage{scalerel}
  \usepackage{lettrine}
  \usepackage{longfbox}
   \usepackage{upgreek}
   \usepackage{stackengine}
   \usepackage{lipsum}
   \usepackage{rotating}
   \usepackage{hhline}

  \definecolor{floralwhite}{rgb}{1.0, 0.98, 0.94}
  
  \usepackage{mathtools}
  \usepackage{accents}
  \usepackage{scalerel}

 \def\vones{\mathbbm{1}}
 
 \makeatletter
 \newsavebox\myboxA
 \newsavebox\myboxB
 \newlength\mylenA
 
 \newcommand*\xoverline[2][0.75]{%
 	\sbox{\myboxA}{$\m@th#2$}%
 	\setbox\myboxB\null
 	\ht\myboxB=\ht\myboxA%
 	\dp\myboxB=\dp\myboxA%
 	\wd\myboxB=#1\wd\myboxA
 	\sbox\myboxB{$\m@th\overline{\copy\myboxB}$}
 	\setlength\mylenA{\the\wd\myboxA}
 	\addtolength\mylenA{-\the\wd\myboxB}%
 	\ifdim\wd\myboxB<\wd\myboxA%
 	\rlap{\hskip 0.5\mylenA\usebox\myboxB}{\usebox\myboxA}%
 	\else
 	\hskip -0.5\mylenA\rlap{\usebox\myboxA}{\hskip 0.5\mylenA\usebox\myboxB}%
 	\fi}
 \makeatother
\usepackage{amsmath,mathtools}

\usepackage[euler]{textgreek}

\usepackage{scalerel,stackengine}
\stackMath
\newcommand\reallywidecheck[1]{%
	\savestack{\tmpbox}{\stretchto{%
			\scaleto{%
				\scalerel*[\widthof{\ensuremath{#1}}]{\kern-.6pt\bigwedge\kern-.6pt}%
				{\rule[-\textheight/2]{1ex}{\textheight}}
			}{\textheight}%
		}{0.5ex}}%
	\stackon[1pt]{#1}{\scalebox{-1}{\tmpbox}}%
}
\DeclareMathAlphabet\mathbfcal{OMS}{cmsy}{b}{n}  

\definecolor{Gray}{gray}{0.85}
\definecolor{LightCyan}{rgb}{0.88,1,1}

\newcommand{\Ex}{\mathbb{E}\hspace{0.05cm}}

\newfont{\sfb}{cmssbx10}
\newfont{\sfbb}{cmr17 scaled 1500}

\def\tran{^{ \mbox{ \scriptsize\hspace{-.22cm} \sf  T\hspace{-.06cm}  }}}

\newcommand\WW{\scaleobj{.75}{\boldsymbol{\cal  W}}}

\def\Aa{\boldsymbol{\cal A}}

\def\PP{\boldsymbol{\cal P}}

\def\LL{ \boldsymbol{\cal L}}

\def\Ccal{ \boldsymbol{\cal C}}

\def\bSigma{{\bf{\Sigma}}}

\def\bLambda{{\boldsymbol{\Lambda}}}

\def\A{\boldsymbol{A}}
\def\B{\boldsymbol{B}}

\def\w{\boldsymbol{w}}

\def\q{\boldsymbol{q}}
\def\x{\boldsymbol{x}}
\def\Pb{\boldsymbol{P}}
\def\Ub{\boldsymbol{U}}

\def\v{\boldsymbol{v}}

\def\u{\boldsymbol{u}}
\def\g{\boldsymbol{g}}

\def\z{\boldsymbol{z}}
\def\f{\boldsymbol{f}}

\def\y{\boldsymbol{y}}

\def\b{\boldsymbol{b}}
\def\I{\boldsymbol{I}}
\def\0{\boldsymbol{0}}
\def\W{\boldsymbol{W}\hspace{-.1cm}}

\def\V{\boldsymbol{V}\hspace{-.03cm}}

\def\D{\boldsymbol{D}}

\def\Q{\boldsymbol{Q}}

\def\L{\boldsymbol{L}}
\def\C{\boldsymbol{C}}

\def\U{\boldsymbol{\cal U}}
\def\Ub{\boldsymbol{U}}
\def\d{\boldsymbol{d}}

\def\bdiag{{\sf bdiag}}

\def\AA{\normalsize {\bf A}}

\def\x{{\bf x}}

\def\ov{\pmb{\mathbbm{1}}}

\newcommand{\be}{\begin{equation}}
\newcommand{\ee}{\end{equation}}

\newcommand{\bq}{\begin{eqnarray}}
\newcommand{\eq}{\end{eqnarray}}
\newcommand{\bqn}{\begin{eqnarray*}}
\newcommand{\eqn}{\end{eqnarray*}}

\newcommand{\ba}{\left[ \begin{array}}
\newcommand{\ea}{\\ \end{array} \right]}

\newcommand{\define}{\stackrel{\Delta}{=}}

\newcommand{\complex}{ \hbox{\rm C\kern-0.45em\rule[.07em]{.02em}{.58em}%
\kern 0.43em}}

\date{}

\usepackage{amsbsy}

\begin{document}

\title{On Distributed Average Consensus Algorithms}

\author{Ricardo Merched\thanks{R.\;Merched  is with the Dept.\,of Electronics and Computer Engineering of  
		Universidade Federal do Rio de Janeiro, Brazil. Email: merched@lps.ufrj.br.
	},  
	\;\emph{Senior Member,~IEEE} 
}

 \maketitle

 \begin{abstract} 

Average\;consensus (AC) strategies play a key role  in every system that employs cooperation by means 
of 
distributed 
computations. To promote consensus,    an  $N$-agent network  can repeatedly  combine    
certain node estimates
until their mean value is reached. Such  algorithms are typically formulated  as  (global) recursive 
  matrix-vector products of size $N$, where  consensus is attained
either   asymptotically or in finite time. We revisit some   existing approaches in these 
directions  and 
propose new 
iterative and exact solutions to the   problem.  Considering directed graphs, this is carried out first by 
interplaying with    
standalone conterparts, while underpinned by the so-called {\em eigenstep} method of finite-time 
 convergence.   Second, we compute the solution  via an exact 
 algorithm that requires, overall,  
as little as 
${\cal 
O}(N)$ additions. For undirected graphs, the latter  compares  favorably to 
existing  schemes that require, in total,   ${\cal O}(K\!N^2)$ multiplications to deliver the AC, where 
$K$ refers to the number of distinct eivenvalues of the underlying graph Laplacian matrix.

\begin{IEEEkeywords}
	Consensus, diffusion, networks, fusion, sparsity, least-squares, adaptation. 
\end{IEEEkeywords}

\end{abstract}

%
%
%
%
%
%




\section{Introduction} \label{sec.finitetimealg}

\IEEEPARstart{L}{east-squares} methods  and their  linearly 
constrained  versions,   are exact in the sense that, 
for a given 
 segment of data, their intrinsic      solution  can be computed in a single shot\,\cite{Ali}.   
 Most 
 of the burden in these computations is   due  to the inversion of covariance 
 matrices estimated from collected  data.  
 In 
 many 
 situations, 
 however, while   data or  covariances  can be fully known,  inversion is prohibited  
   owing to complexity limitations  and/or because of structural constraints 
 on the inverse itself. A classical example of such difficulties arises in a scenario of spatially 
 distributed parameters, which includes  
 {\em consensus} and {\em diffusion} as popular learning     strategies\,\cite{Ali1}. Such schemes are 
 essencially {\em linearly 
 constrained}  algorithms, as in view of \cite{Frost}, albeit formulated in terms of a 
 global 
 network 
 model\,\cite{mer22A},\cite{mer22B}.  While the respective constraining matrix\,($\C$)  in  these  
 formulations   may or may not exhibit a sparse structure\,(in accordance with the topology of the 
 associated   
graph),  should the   solution be expressed in terms of its  
pseudoinverse ($\C^{\dagger}$), it will generally exhibit filled-in entries, and, as a result,  
  introduce new  
challenges  for its proper, practical realization.

To overcome  such complexity and structural limitations, one can rely solely on  gradient-based 
updates that do not require inversions. For instance,  an {\em average\,consensus} (AC) becomes quite 
useful 
when   
       interleaved with self-learning algorithms in cooperative networks  --- see, e.g., 
       \cite{Ali1},\cite{mer22B}. To avoid inversions, the author in the latter makes
use of the well-known
conjugate-gradient\,(CG)\,\cite{Gill} method,  as a means to indirectly realize the AC step 
of the   solution. The advantage is  that,    given  a set of values collected from an 
$N$-node  
network,   this could  in theory    be  
       reached   in
       at most   $N$  	iterations,   were it not for the dependency    
    of  an optimized, data-dependent
       global parameter.

       A fast computation of the AC is of great interest, since it is inherent to   distributed 
       processing in 
       several 
       other 
       applications. The subject has been long investigated, as e.g., in    
         \cite{Xiao},\cite{Boyd}, 
         considering optimized time-invariant  matrices,    possibly  bearing negative  combination 
         weights.    For 
         example, 
         the CG 
         method   was 
         also used in  
\cite{PP1} in order to replace  covariance  inversions by linear operations in the processing of signals 
over     
distributed   arrays. In particular, in addition to CG,  the authors in \cite{PP1}  rely heavily on   
a so-called {\em finite-time  AC} algorithm derived  in \cite{Moura2014}, where it reemerged  as a 
soubroutine for 
the  
computation of several 
averages that appear as parts of their distributed computations. The same idea  
was also employed in 
 \cite{Kibangou} in the scenario of  sensor networks,  yet assuming non negative combination coefficients.   
The notion of finite-time 
       convergence   dates back to 
       \cite{Ryan},\cite{Bhat}, where stability conditions   
       were   studied  for continuous-time  autonomous systems --- see also 
       \cite{Cortes},\cite{LWang},\cite{Hendrickx}. The concept of consensus strategies, on 
       the other 
       hand,  was ushered in much 
       earlier by Degroot\,(1974)
       in \cite{Degroot}, with further roots appearing in  the  references therein.

It turns out that the finite-time AC algorithm in \cite{Moura2014} was not new. It  had been 
   proposed  a decade earlier by Battaglia, in \cite{Batt} (within his master's thesis), as a special case 
   of  
   the well-known
   standalone,    
 {\em vanilla}  
steepest-descent algorithm. 
 The recursion was referred to as the {\em eigenstep method}, and can be simply  recognized  as a sequence 
 of   
  gradient-based iterations,   with   time-varying 
 step sizes   selected particularly   as the reciprocals of the    underlying  
 Hessian matrix's  eigenvalues.	  As a 
 result, in the context 
 of graphs, the  algorithm takes the form of at most $N\!-\!1$  products  of  
 distinct   matrices at the network level (since, in  the setting of \cite{Moura2014},  these  were 
 restricted  
 to 
 Laplacian 
 structures  exhibiting  a simple zero eigenvalue).  
  We remark that in\;\cite{Safavi},  the authors also claim novelty on this very same idea,  by  referring to 
  it    as a  method of 
    ``{\em successive nulling of 
    eigenvalues}''. This is however,  the same  eigenstep method introduced in  \cite{Batt},  
    though formulated  in the   distributed  sense.

Motivated by the need of efficient AC schemes,  we pursue in this work a few  
 strategies that aim  improved convergence and low
complexity, for algorithms derived from  iterative and exact 
formulas.  The  ideas   hereinafter  are  organized as follows: 

\vspace{-.2cm}
 \begin{enumerate} 
 \item  ({\em Sec.\,\ref{backgr}}) To ensure that this presentation is  self-contained, we 
 revisit the AC concept 
 in a 
 more   unified manner, by 
  introducing it as a particular application of a simple  linearly-constrained problem. We 
highlight the fact that its solution expresses the AC exactly from   values stored  by the   nodes 
 of 
 a   graph. In this sense, one can borrow from an extensive literature of  iterative algorithms, 
 possibily with the 
 use of momentum terms\,(e.g.,  \cite{Nesterov},\cite{Qian})  in order to 
 achieve AC  efficiently in a distributed manner. This includes the eigenstep method \cite{Batt}, for which   
 we point out the
 equivalence to the algorithms in  \cite{Moura2014} and \cite{Safavi}, while  extending this notion 
 to directed networks;     \vspace{.1cm}
 
%

 \item  ({\em Sec.\,\ref{ACalg}\;$\&$\;\ref{discuss}}) We argue that, in general, the eigenstep iterations 
 can 
 become quickly unstable, specially for large and 
 ill-conditioned matrices, in which case  the AC cannot be met. We address this issue by exercising some 
 natural  
 interplays between   standalone and 
 distributed   schemes,  in light of their intrinsic 
 power iterations. This allows us to implement the eigenstep method in a rather stable manner,    
  shedding  
  additional light into the construction of  finite-time convergence  
 algorithms.  An important conclusion from this discussion is that, without  constraining the weights to some 
 form of normalization,  
 convergence of the power iterations  to AC can occur  faster compared to a scenario where the
 Laplacian coefficients are  optimized\,\cite{Xiao},\cite{Boyd}, as verified through simulations. In fact, 
 even a  Nesterov accelerated 
         gradient\,(NAG)  method\,\cite{Nesterov} applied to problems under a non-optimized Laplacian 
   can    outperform  one with optimal  weights, considering an acceptable error margin\,(after,  
         typically,  
         $2N$
         iterations); \vspace{.1cm}

 \item   ({\em Sec.\,\ref{sec.exact}}) With   focus  on reduced complexity and stable solutions, we show, for 
 undirected graphs, 
 that the exact AC  can be 
 obtained 
 efficiently in $2N$ iterations, via   
 back/forward-substitution recursions at the network level. As a byproduct, the exchange of information can 
 be minimized, and  we can attain an  
 overall 
 complexity    as low as ${\cal O} (2N)$ additions. This can  be contrasted with the eigenstep method, which 
 would require ${\cal O} (K\!N^2)$ multiplications in total\,(where 
 $K$
 is  the number of distinct eivenvalues of the underlying graph Laplacian), or a power iteration, which may 
 take too long to 
 converge depending on the conectivity of the graph. From this same development, we obtain 
 a second algorithm that achieves AC  in the same number of iterations as the eigensteps, 
 and requiring, overall, a complexity of  ${\cal O} (N^2)$ additions.    
  \end{enumerate}
 
  While this presentation deals with independent AC iterations, we highlight that its use is also called for 
 upon its
  combination
  with 
  self-learning recursions, as  it  may be key to minimizing the 
  complexity of the AC 
  step ---- see {\em Sec.\,\ref{ACdiffu}}.

  \noindent {\bf Notation\,}:
    	 We  employ regular, italic, and caligraphic  types, at times for the same letter, e.g., $\x$, 
    	 $\boldsymbol{x}$, or   $\Aa$, $\A$, $\AA$, ${ \sf A}$. Vector entries are denoted as 
    	 $\x(k)$,  and  submatrices by  colon notation, as e.g., in   
    	 $\A(m\!:\!n,p\!:\!q)$. 
    	 The complex conjugate transposition is written as $\A^{\star}$, while $\A^\dagger$ denotes the 
    	 pseudoinverse.  
    	 We also define the pinning vector ${\sf e}_k=[\,0 \;\, 0 \,\cdots \,0\,\;1  \,\;0\,\cdots 
    	 0\,]\tran$ with unity  
    	 $k$-th 
    	 entry, the all-ones vector  $\ov =\left[\,1\;\;1\;\;\cdots\;\;1\,\right]^{\sf 
    	 T}$   
    	 of appropriate dimensions, and a ``block'' of zeros  
    	by $\0$.
    	  The operator $\bdiag\{\cdot\}$  forms a block diagonal matrix from a 
    	 block-column vector,  extending the notation $\text{diag}\{\cdot\}$, 
    	 used for scalars.  The operator  $ \text{lower}\{{\sf A}\}$  extracts the   lower triangular 
    	 part of ${\sf A}$.   

 
 \section{Background and General Setting} \label{backgr}
 
Assume we  know a  given   set of  numbers   collected  into an ($N\times 1$) vector  at time $i=0$,  
say,  
 $\w_0 \in \complex^N$, and 
 consider the following   minimum-(Mahalanobis)\,norm problem:
  \be
    \min_{w} \;\; \left\| \w  -  \w_0 \right\|^2_{{\Pb}^{-1}_{\!0}}   \;\;\;\;\text{subject to}\;\;  
    \C'\!\w  =    \d  \label{compfcccwo}
   \ee
 where $\C'$ is $(N\!-\!R) \times N$. The ($N \times N$) covariance matrix ${\Pb}^{-1}_0$ measures the  
 uncertainty and  
 correlation 
 among the 
 entries of $\w_0$, which is subject to a 
   linear constraint.  Several stochastic and deterministic constrained objective functions can be cast into 
   this 
 form  via {\em completion-of-squares}, which include  least-mean-squares\,(LMS),  
 Affine Projection, and 
 deterministic least-squares problems\,\cite{Ali}, all leading  to the same form of a  (weighted) 
 minimum-norm 
 solution, denoted by $\w^o$:
  \be 
     \w^o\;=\;  \w_0 +     \Pb_{\!0}\C^{'\sf T}\Big( \C' \Pb_{\!0} \C^{'\sf T}   
    \Big)^{\!-1}    
   \!\!( \d - \C' \w_0 )    \label{brlsmnl000i}
   \ee
   Depending on the cost at hand, the ``uncertainty'' $\Pb_{\!0}$ can be  updated  from a previous 
   estimate, 
   say, 
   $\Pb_{\!-1}$, 
   corresponding to some prior guess,  
   $\w_{-1}$. 
 If  $\Pb_0$ is not available,   we can set $\Pb_0=\I$\;[see, e.g., (\ref{brlccccnlb})  in 
   Sec.\,\ref{ACdiffu}]; otherwise, we can incorporate $\Pb_0^{1/2}$ into  
 $\C'$ 
 in 
 (\ref{compfcccwo}) by defining $\C \define \C' \!\Pb_0^{1/2}$, and replace the corresponding parameter 
 by  
 $\Pb_0^{-1/2}\w \longrightarrow \w$.  Moreover, by  denoting  $\f \define    \C\tran ( \C   
 \C\tran 
 )^{\!-1}\!\d =  \C^{\dagger} \d $, the 
 formula
 (\ref{brlsmnl000i}) can be 
 rearranged as  
  \bq
    \w^o   &\!\!\! =&  \!\!\! \Big[\I -    \C\tran \!\left( \C   \C\tran 
     \right)^{\!-1}  \!\C\Big] \w_0   \,+\, \f    \label{brlsmnl0002}  
    \eq
  Expression (\ref{brlsmnl0002}) is simply the general solution 
 of 
 a  
 linear 
 system 
 of equations $\C\w=\d$ in terms of  the 
  orthogonal projection of a 
 particular vector $\w_0$  onto the kernel of  $\C$. 
 It is also instructive to realize that (\ref{brlsmnl0002}) 
 can be 
 seen as a single 
 update
 of a {\em block normalized} LMS algorithm,  from the prior guess,   $\w_0$.  
 
 Now, consider the singular-value-decomposition\,(SVD),
 \be 
 \C = \Ub \ba{cc}\bSigma & \0\ea \Q\tran  \label{SVDofC}
 \ee
 It is well known that the projection in (\ref{brlsmnl0002}) can be written as 
 \be 
 \I -    \C\tran \!( \C   \C\tran )^{\!-1}  \!\C = \Q \ba{cccc} \!\!\!0 &&& \\ &  \!\!0 &&\\ && \ddots & \\ & 
 &&  \!\I 
    \!\!\!\ea 
    \Q\tran\,=\,\V_{\!N}\V_{\!N}\tran \label{DecQDQs}
 \ee
 in terms of the   eigenvectors matrix $\Q\define [\v_1\;\;\v_2\;\cdots\;\v_N]$, or,  via its
 submatrix   
 $\V_{\!N}$ 
 comprising  its right-most $R$  eigenvectors.
 That is, in general,   (\ref{brlsmnl0002}) can be expressed as 
 \be 
 \w^o= \V_{\!N}\!\;\!\V_{\!N}^{\sf T} \w_0\,+\, \f \label{genexpressfvdvv}
 \ee
  For example, 
  if 
  $\C$ 
  is an $(N-1) \times N$, full row rank matrix, then $\V_{\!N} = \v_N$ becomes  a column vector. By 
  associating 
  $\v_N$ 
  to a certain 
  probability distribution $\{p_k\}$, $k=1,\ldots,N$,  each entry $[\w^o(k)-\f(k)]/p_k$ will  correspond  to  
 the
  mean value of $\w_0$, denoted by 
  $\bar{w}_0=\v_N\tran \w_0$. For instance, in the context of cooperative networks, we can associate  $\w_0$ 
  to an 
  extended 
  vector that contains  estimates  read from all of its $N$  nodes\,(see Fig.\,\ref{fig.net_asymm_20}, 
  further 
  ahead). 
  In 
  the 
  case of adaptation and diffusion   
  strategies,  
  nodes  
  first 
  self-learn 
  their parameters, so that   
  the exact linearly-constrained solution  is obtained by     averaging  all  these estimates  
  across 
  the  
  network\;[i.e., in this scenario, these strategies  simply amount to    constrained  adaptive
  algorithms in the linear sense ---- see these steps, e.g., as  (\ref{brlccccnlbx}),(\ref{brlccccnlb}), 
  Sec.\,\ref{ACdiffu}].
  
  The structure of the  projection (\ref{DecQDQs}) leading to the solution in (\ref{genexpressfvdvv}) 
  suggests 
  that  a 
  more general, asymmetric  rank-one mapping from $\w_0$ to $\w^o$, should also hold. Indeed, consider two 
  distinct full row rank matrices $\C_1$ and $\C_2$. Then,  
   \be 
    \underbrace{\I -    \C_2\tran  ( \C_1   \C_2\tran )^{\!-1}  \C_1}_{\A}  \,=\,\V_{\!N}\Ub_{\!N}\tran 
    \label{DecQDQsgen}
   \ee
   so that (\ref{genexpressfvdvv}) is  replaced by 
    \be 
    \w^o= \V_{\!N}\!\;\!\Ub_{\!N}^{\sf T} \w_0\,+\, \f \label{genexpressfvdvv1}
    \ee

\noindent {\em Proof}\,: Without loss of generality,  assume that $\{\C_1,\C_2\}$ are  $(N-1) \times N$ 
matrices, since this will be the case of interest here.  This is required, for example, if we are to 
   linearly   interelate 
    all the 
    elements of  $\w$. For example, in the distributed scenario, this reflects the fact that the 
    network is 
    {\em 
    connected}.    Next, consider the following partitions:   
    \be 
   \C_1\,=\,\ba{cc} \xoverline{\C}_1 &  \b\tran \ea \,,\;\;\;\; \C_2\,=\,\ba{cc} 
   \xoverline{\C}_2 
  &  \g\tran \ea \label{Cpartix}
   \ee
 The mapping $\A$ indicated in (\ref{DecQDQsgen}) can be simplified as follows:
        \bq
         \A&\!\!\! \!  =& \!\!\!\! \I -    \ba{c} \!\! \xoverline{\C}\tran_{\!2}  \\   \g \ea\!\left(   
          \xoverline{\C}_1  \xoverline{\C}\tran_2 + 
         \b\tran 
         \g\right)^{\!-1}    
        \!\ba{cc}  \xoverline{\C}_1 &   \b\tran \!\ea   \label{exctsolLL11} \\
        &\!\!\!\!  & \!\!\!\!\I -    \left[\!\begin{array}{cc} \!\! \xoverline{\C}\tran_2  \\   \b  
        \end{array}\!\!\right]\!\!\left(  
        \!\xoverline{\C}_2^{-\!{\sf T}}\xoverline{\C}_1^{-1} 
        \!-\! 
        \dfrac{\xoverline{\C}_2^{-\!{\sf T}}\!\xoverline{\C}_1^{-1} \b\tran 
        \!\g \xoverline{\C}_2^{-\!{\sf T}}\xoverline{\C}_1^{-1}}{1+\g 
        \xoverline{\C}_2^{-{\sf T}}\xoverline{\C}_1^{-1}\b\tran}\!\right)  
        \! \!\left[\begin{array}{c} \!\!\xoverline{\C}_1\tran  \\  \!\!\b \!\end{array}\!\!\right]^{\!\!{\sf 
        T}}   \nonumber \\
        &\!\!\!\! =&\!\!\!\! \dfrac{1}{1+\g 
                \xoverline{\C}_2^{-{\sf T}}\xoverline{\C}_1^{-1}\b\tran} \ba{cc} 
        \xoverline{\C}_1^{-1}\b\tran \! \g \xoverline{\C}_2^{-{\sf T}}   & - 
        \xoverline{\C}_1^{-1}\b\tran    \\  - \g\xoverline{\C}_2^{-{\sf T}}    &  1 \ea  \nonumber  \\
        &\!\!\!\! =&\!\!\!\! \gamma \ba{c}  \! \xoverline{\C}_1^{-1}\b\tran \hspace{-.2cm}   \\ -1 \ea 
        \ba{cc} 
        \g \xoverline{\C}_2^{-{\sf T}}  
        & 
        -1 \ea 
        \eq
        where we further  defined the scalar factor
        \be 
        \gamma \,\define\,(1+\g 
                        \xoverline{\C}_2^{-{\sf T}}\xoverline{\C}_1^{-1}\b\tran)^{-1}
        \ee
        
        As a result, the extended vector in (\ref{genexpressfvdvv}) becomes
        \bq
        \!\!\!\!\!\!\!\!\w^o &\!\!\!\! =&\!\!\!\! \gamma \ba{c}  \! \xoverline{\C}_1^{-1}\b\tran 
        \hspace{-.2cm}   
        \\ 
        -1 
        \ea 
                \left[\!\!\begin{array}{cc} 
                \g \xoverline{\C}_2^{-{\sf T}}  
                & 
                -1 \end{array}\!\!\right]
        \w_0 \;+\;\f
        \label{expswow0} \\
         &\!\!\!\! =&\!\!\!\! \gamma \underbrace{\left[\!\begin{array}{cc} \!\! \xoverline{\C}_1^{-1} & \0 \\ 
         \0 & 1 
        \end{array}\!\!\right]}_{\W_1^{\;-1}} 
        \!\underbrace{\left[\!\begin{array}{cc} \!\! \b\tran \\   -1   \end{array}\!\!\right]}_{\b_1\tran}   
        \!\underbrace{\left[\!\begin{array}{cc} \!\!\g & -1 
          \end{array}\!\!\right]}_{\g_1}  
        \!\underbrace{  \left[\!\begin{array}{cc} \!\! 
       \xoverline{\C}_2^{-{\sf T}}  & \0 \\ \0 &  1 \end{array}\!\!\right]}_{\W_2^{\;-{\sf T}}}  \w_0 +\f 
       \nonumber 
        \eq
       yielding, in terms of the indicated variables above,
        \be 
         \w^o \,=\, \gamma  \W^{\,-1}_1 \b_1\tran \g_1 \W_2^{\,-{\sf T}}   \w_0 \; +\;\f \label{prodmatswog}
        \ee
       Comparing this expression with (\ref{genexpressfvdvv1}), we  identify $\{\V_{\!N},\Ub_{\!N}\}$ as the 
       column 
       vectors 
       \be 
       \v_N \;=\; \sqrt{\gamma}\W^{\,-1}_1 \b_1\tran\;\;\;\text{and}\;\;\; \u_N \;=\; 
       \sqrt{\gamma} \W^{\,-1}_2 
       \g_1\tran
   \ee

  Now, inspired by (\ref{brlsmnl0002}), let us 
 drop 
 the corresponding  ``normalization'' term $( \C_{\!2}  
 \C_1\tran )^{\!-1}$ in (\ref{DecQDQsgen}), and replace $\w^o$  by an estimate, 
 $\w_i$, 
    obtained recursively via  a time-varying step-size, $\mu_i$, as
  \be
    \w_i \,=\,  \w_{i-1} -    \mu_i \C_{\!2}\tran    \!\C_1 \w_{i-1}   \,+\, \f_i  
    \label{brlsmnl0blms}  
   \ee
 where $\f_i= \mu_i \C_{\!2}\tran \d$. This is commonly motivated in case $\C_1=\C_2=\C$, as an 
 unnormalized   gradient 
 based  
 algorithm. Here, we  consider the above extension   
 and  seek the 
 values of $\mu_i$ that will  allow us 
   to achieve  the solution $\w^o$ asymptotically, or possibly  in a finite number of iterations. To this 
   end, it is a  
   standard procedure   to 
   define the error quantity  $ \widetilde{\w}_i = \w^o- \w_i $, 
 so that substituting it into (\ref{brlsmnl0blms}) we obtain 
  \be
 \widetilde{\w}_i \,=\, \widetilde{\w}_{i-1} \,+\, \mu_i \C_{\!2}\tran    \!\C_1  (\w^o 
 -\widetilde{\w}_{i-1}  
 )   
 \,-\, \f_i   
 \label{brlsmnl0blmsm}  
  \ee
 or, equivalently,
  \be
 \widetilde{\w}_i \,=\, \Big(\I -\mu_i \C_{\!2}\tran    \!\C_1  \Big) \widetilde{\w}_{i-1}    + \mu_i 
 \C_{\!2}\tran\!\C_1 \w^o    
 \,-\, \mu_i  \C_{\!2}\tran \d  
 \label{brlsmnl0blmsm1}  
  \ee
 Since $ \mu_i \C_{\!2}\tran    \!\C_1 \w^o  -\mu_i \C_{\!2}\tran \d   =\mu_i \C_{\!2}\tran(   \C_1\w^o  -  
 \d)=\0$, 
 we have that 
    \be
 \widetilde{\w}_i \,=\, \left(\I -\mu_i  \C_{\!2}\tran\!\C_1\right) \widetilde{\w}_{i-1}       
 \label{brlsmnl0blm2}  
  \ee
whose solution  at an arbitrary instant $i$ is given by
 \be 
 \widetilde{\w}_i \,=\, \prod_{k=1}^i \left(\I-\mu_k  \C_{\!2}\tran\!\C_1\right) \widetilde{\w}_0 
 \label{sfeteruis}
 \ee
in terms of the initial error vector
 \be 
 \widetilde{\w}_0 \,=\, \w^o - \w_0 \,=\,  \C_{\!2}\tran \left( \C_1   \C_{\!2}\tran \right)^{\!-1} \! ( \d - 
 \C_1 \w_0 )    
 \label{wtildeCtC}
 \ee
 
In general, we have that $\widetilde{\w}_i \rightarrow \0$ (and bounded away from 1),  provided that the 
eigenvalues 
of the product in 
 (\ref{sfeteruis}) are stable. Moreover, note that at each iteration, we can attempt to set one of its 
 eigenvalues 
 exactly to 
 zero by choosing the step size $\mu_i$ accordingly.  That is, similarly to what is implied by 
 (\ref{SVDofC}) when $\C_1=\C_2$, assume that 
   a more general  diagonalization of the form
    \be 
    \C_{\!2}\tran \!\C_1 \,=\,   \Q \bLambda \Q^{-1}
    \label{SVDlamsrd}
    \ee 
exists,  such that
 \be 
 \bLambda=\text{bldiag}\{\lambda_1\I_{\!m_1},\lambda_2\I_{\!m_2},\ldots,\lambda_{K}\I_{\!m_{\!K}},0\}  
 \label{lmultipl}
 \ee
 In this case, the matrix $\C_{\!2}\tran \!\C_1 $ will exhibit  one eigenvector  
 associated to the 
 eigenvalue   0, and $K$ eigenvalues with multiplicities $m_1,m_2,\ldots, 
 m_{K}$,  respectively. That would happen in a distributed scenario, e.g., when   the network is connected. 
 Now, assume that $\{\lambda_k\} \neq 0$, and set  $1- \mu_k \lambda_k =0$, for $k=1,2,\ldots,K$. 
This means that we can  reach the exact solution $\w^o$ in $K$ iterations of (\ref{brlsmnl0blms}), by 
selecting 
 \be 
 \mu_i\,=\, \dfrac{1}{\lambda_i}\,,\;\;\;\;\text{for }i=1,2, \ldots,K. \label{muk1ovlam}
 \ee
 
 \noindent {\em Proof}\,: Using (\ref{muk1ovlam}), we express (\ref{sfeteruis}) at time 
 $K$ as
 \bq 
 \widetilde{\w}_{K} &=&  \Q\left(\prod_{k=1}^{K} \left(\I-\mu_k \bLambda \right)\right) \Q^{-1}
 \widetilde{\w}_0 
 \label{QproQtra} \\
 &=& \Q \ba{cccc} \0 &&& \\ &  \0 &&\\ && \ddots & \\ & && 1 \ea \Q^{-1} \widetilde{\w}_0   \nonumber \\
 &=& \v_N\u_N\tran\widetilde{\w}_0   
 \;\;\;\;\;\;\;\;\;\;\;\;\;\;\;\;\;\;\;\;\;\;\;\;\;\;\;\;\;\;\phantom{..} \text{[using
  (\ref{genexpressfvdvv})]} \nonumber  
  \eq 
  where 
  \be 
  \u_N \define  {\sf e}\tran_N\Q^{-1}  \label{uNQinv}
  \ee 
   in terms of the pinning vector $ {\sf e}_N$. Then, using 
   (\ref{wtildeCtC}), we obtain
 \be 
 \widetilde{\w}_{K} \;=\;\v_N\!\underbrace{\u_N\tran  \C_{\!2}\tran}_{=\0} \left( \C_1   \C_{\!2}\tran 
 \right)^{\!-1}  \!( 
 \d - 
  \C_1 \w_0 ) \;=\; \0 \nonumber 
 \ee
 The latter is justified by multiplying (\ref{SVDlamsrd})   from the left by  $\u_N$, which gives $(\u_N 
 \C_{\!2})\tran\C_1=\0$. Since $\C_1$ is assumed full row rank, this further implies that $\u_N 
 \C_{\!2}\tran=\0$, 
 yielding $\w_{K}=\w^o$.  
 
 When $\C_1=\C_2=\C$, the corresponding gradient-descent recursion that makes use of the  time-varying step 
 size (\ref{muk1ovlam}) is referred 
 to 
 as  
 the {\em 
 eigenstep} method. The idea  was introduced for the first time  in \cite{Batt}, for   
 general-purpose, 
 possibly full 
 rank,   Hermitian
 matrices\,(i.e., replacing $\C\tran \!\C$ above, in which case the algorithm would 
 terminate in at 
 most $N$ 
 iterations). 
 
\section{Fast Average Consensus\,(AC) Algorithms} \label{ACalg}

\begin{figure}[htb]
   			      \epsfxsize 15.4cm \epsfclipon
   			      	\begin{center}
\includegraphics[clip,width=5cm]{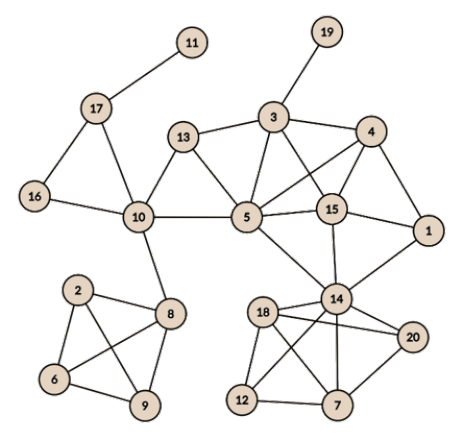} 
   			      		  \caption{Cooperative network wit  $N=20$ nodes, arbitrarily 
   			      		  numbered.}\label{fig.net_asymm_20}
   			      	\end{center}
   			      \end{figure}

   	Algorithms that  converge  in finite-time  are  relevant in the realm of  cooperative networks, where  
   	the 
   			constraining matrix  $\C$ is known,  and often chosen (although not necessarily) with the 
   			sparsity 
   			pattern that reflects  a given network topology. When this is the case,  the  rows 
   			of 
   			$\C$ comprise  
   			certain weights that can be  used to combine   estimates  from  
   			neighboring agents, which in turn can be described by an extended vector, $\w$,  having
   			vector-valued entries. This is 
   			illustrated in Fig.\,\ref{fig.net_asymm_20}, for an undirected  
   			network of agents  comprising 
   			$N=20$ nodes. 
   			 
   			   			\vspace{-.3cm}
  
  For simplicity of exposition, let $\w_0$ be a vector of {\em scalar}-valued coefficients      
  numbered in an arbitrary way. The 
  objective 
  in   AC  strategies is to deliver, at 
    each node,  the arithmetic mean of all coefficients,  $\bar{w}_0 
  =(\sfrac{1\!}{N})\vones\!\tran 
  \w_0$,    
  according to some recursive algorithm. That is,  the goal is  to compute
  \be 
  \w^o\,=\,  \frac{1}{N} \vones \!\vones\!\tran \w_0 \,=\,\bar{w}_0\vones  \label{averwow0}
  \ee
  
   In many AC strategies found in the literature, the algorithm
  entails  network-level iterations  of the linear form  
   \be 
   \w_i\,=\, \A_i \w_{i-1} ,\;\;\;i=1,2,\ldots,N_o \label{tviter}
   \ee 
   which can result, e.g., from  a  
  linear approximation problem of designing   a set of  (time-varying) matrices $\{\A_i\}$, such that 
  \be 
  \w^o\,=\,\left(\prod_{i=1}^{N_o} \A_i \right) \!\w_0  \label{prodsetW}
  \ee
  To achieve consensus,  the choice of $\A_i$ 
      must be   such that 
      \be 
      \left[\prod_{k=1}^{j_o}\A_i \right]_{\ell k} \!\neq \;0\;\;\;\text{for some }j_o>0 \label{Alkcond}
      \ee
     otherwise,   the above product   will not  effectively operate on all entries of 
     $\w_0$. Thus, in 
          order to reach $\bar{w}_0$ with  the smallest $j_o$ possible, the nodes need to store a set of  
          optimized, possibly time-varying coefficients, according to some criterion.   
     In other words, for a connected network, one  
     would 
     like 
     to 
     design  sparse   matrices  
     $\A_i$ that reflect  the adjacency  structure,   while satisfying  (\ref{Alkcond})  considering  a 
     minimum  
     computational burden.
%
%
%
%

  One  route in this direction is to first enforce the structure of $\w^o$ in  
  (\ref{averwow0}) by considering a linear constraint. That is, let us explicitly denote the vector   
  $\vones$ with its size, i.e., $\vones_{\!N}$. Then, given that  for any $(N-1) \times N$ matrix 
  $\B$, we have
  \be 
  \B \w^o= \bar{w}_0  \B\vones_{\!N}   = \bar{w}_0 \b = \bar{w}_0 \text{diag}(\b)\vones_{\!N-1}     = 
  \left[\text{diag}(\b)\;\;\0\right]\!\w^o
  \ee 
  or, equivalently, that $\big(\B-\left[\text{diag}(\b)\;\;\0\right]\big)\w^o=\0$, by defining 
  \be 
  \C \define  \B-\left[\text{diag}(\b)\;\;\0\right] 
  \ee
      we can formulate  (\ref{compfcccwo}) with $\d=\0$, i.e., with $\C 
      \w=\0$. Since  the entries of $\w$   must be identical, it further 
  implies that $\C \vones = \0$. In a connected network, $\C$ is $(N-1) \times N$, so that multiplying this 
  relation 
  from 
  the left by $\C\tran$\!, we further conclude that this is  equivalent to
  selecting 
   a matrix $ \C\tran \!\C$ having  a (scaled) eigenvector $\vones$  associated to the simple 
  eigenvalue 
  0. As a result, the mapping in   (\ref{genexpressfvdvv}) is satisfied with    
  $\v_N=(\sfrac{1\!}{\!\sqrt{\!N}})\vones$, which can be achieved in several 
  ways.    One possibility is to set $\C\tran \!\C = \LL$,   a  Laplacian 
    matrix having the general form
    \be
              [ {\LL}]_{k\ell} \;=\;\left\{ \begin{array}{ll} \!\!\sum\limits_{\ell \neq k} c_{k\ell}, & 
              \text{if 
              }k=\ell 
              \\\!\!- c_{k\ell},  & 
              \text{if 
              }k\!\neq\!\ell    \text{ with nodes $k$ and $\ell$ $\in {\cal N}_k$}
               \\ \!\!0, & \text{otherwise} \end{array}   
              \right. 
              \label{LLaplungen}
              \ee
    which satisfies
  \be 
  {\LL}\vones\,=\,\0 \label{Lii0}
  \ee 
    For example, the said {\em unormalized}
    Laplacian with constant weights is defined as
    \be
        [ {\LL}]_{k\ell} \;=\;\left\{ \begin{array}{ll} \!\!n_k\!-\!1, & \text{if }k=\ell \\\!\!-1,  & 
        \text{if 
        }k\!\neq\!\ell  
        \text{ with nodes $k$ and $\ell$ $\in {\cal N}_k$}   
         \end{array}    \right.  
        \label{LLaplun}
        \ee 
       otherwise, $ [ {\LL}]_{k\ell}=0$, with   $n_k$  denoting    the degree of node $k$.
  
  Now, since   $(\C \C\tran)^{-1}   $ in (\ref{brlsmnl0002}) is a fully dense matrix,  in general
   its structure 
      does not 
      conform with the graph topology, and   a solution  of the form (\ref{brlsmnl0002}) must be 
      realized  
     according to some suitable strategy.   One possibility is to resort to the algorithm 
     (\ref{brlsmnl0blms}), for which  
      the  network iterations assume the desired form  
     (\ref{tviter}), with   $\A_i=\I-\mu_i{\LL}= \I-{\LL}_i $. Moreover, for a fixed matrix ${\LL}_i 
     =\LL$,  we can seek an optimal, time-invariant, Laplacian matrix ${\LL}^o$, 
         that yields the fastest distributed  averaging,   by minimizing the spectral norm   $\|\I - 
         \LL\|_2$.  This 
                  would require less information to be stored by the nodes with guaranteed stability, 
                  although under  
                  asymptotic 
                  convergence.  That is, 
         in view of (\ref{Lii0}), and considering the inherent sparsity pattern $\mathbb{S}$    of ${\LL}$,  
         we 
         can 
         solve 
             \be 
             \min_{{\LL}} \left\|\I - {\LL}\right\|_2\;\;\;\;\text{s.\,t. }\;  {\LL}\in \mathbb{S},\;   
             {\LL}={\LL}\tran,  \; {\LL}\vones = \0  \label{probAfor}
             \ee
        which 
         reduces to the minimization of  
            $\text{max}\{1\!-\!\lambda_2,\lambda_N\!-\!1\}$ ---- see \cite{Boyd},\cite{Boyd2}.  The cost 
            (\ref{probAfor}) consists in a convex optimization, and was
            formulated in the latter references straightforwardly as 
                       a standard semidefinite programming,   along with other useful related objective 
           functions.\!\footnote{For   fixed  $\A_i =\I- \mu{\LL}$, we 
           have 
           that $\mu 
            < \mu_{\sf max} =2/\lambda_{\sf max}({\LL})$, while its optimal value collapses to   $\mu^o= 
            2/(\lambda_2+\lambda_N)$.  In this case, convergence  depends on the eigenvalue  spread of $\LL$, 
            which  can  be extremely 
                      slow.}

  As an alternative route, one can  resort to the eigenstep procedure of Sec.\,\ref{backgr}, by employing  a 
  time-varying matrix 
   $\A_i$ with (\ref{muk1ovlam})  in   (\ref{tviter}), so that  
   (\ref{averwow0}) is computed    in  exactly
   $K$ iterations of the algorithm. That is, we define  
   \be 
   \A_i\,=\, \I -   \dfrac{1}{\lambda_i}{\LL} \;\;\;\;i=1,2,\ldots,K \label{weightsAil}
   \ee 
 so that $\w^o= \w_{K}$. 
 This    is attractive in the   sense that,   compared to  an exact AC expression that would in 
 principle require $N$ iterations of a sequential 
   algorithm\,(without any structural constraint), or  a slow-convergent iterative scheme, the same 
   solution can be reached faster,  
by relying on  
   the 
      possible eigenvalues 
      multiplicity in ${\C\tran \!\C}$.

     Still, there are some caveats related to this approach.  First, even though the algorithm aims to 
     achieve convergence 
   exactly, it is  prone to become numerically unstable  depending on the size and conditioning of 
   $\C$.  Note that while the modes of $\A$ remain bounded for fixed step sizes,   the eigensteps  
   employ  
    values of $\mu_i$  that can exceed $\mu_{\sf max}$  by far, in the hope that at the final 
   iteration     
   there will be no   accumulation of numerical errors\,(i.e., by assuming all modes in (\ref{QproQtra}) are 
   set 
   precisely to zero). This,  
   however, is not the case in 
   general, and the fact that convergence can be achieved   exactly at the $K$-th iteration does not mean 
   that, 
   before that,   we will be close to the solution. The latter may be relevant, for example, in a scenario 
   of ``adaptation-and-diffusion'', where  only a single pass of combinations  in ${\cal N}_k$   is allowed  
   after 
   agents self-learn 
   their estimates\,(see  Sec.\,\ref{ACdiffu}, for the case of LMS adaptation). 
   
   To see this, suppose that upon defining   the step sizes in (\ref{muk1ovlam}) we incur a negligible, yet 
   nonzero 
   error, such that the corresponding modes with respect to the error  vector are also nonzero, i.e., $ 
   1-\mu_k\lambda_k\approx \epsilon_k$.  
   Assume  for simplicity that $N\!=\!5$, and that all eigenvalues are distinct. From  (\ref{QproQtra}), 
   we obtain
      \bq
        \prod_{k=1}^{N-1} \left(\I-\mu_k \bLambda \right) & \!\!=&  {\small \!\!\!\!\ba{ccccc} \!\epsilon_1 
       & 
       & 
       & 
       &\\ &\!\! \!\!\!(1\!-\!\mu_1\lambda_2)  & & & \\ & 
      & 
      \!\!\!\!\!\!(1\!-\!\mu_1\lambda_3)& & \\ 
      & &&    
      \!\!\!\!\!\!(1\!-\!\mu_1\lambda_4) &  \\ & & & &\!\!\!\!\!\!\!\!\!\!\!1 \!\!\!\!\!\!\!\!\!\ea}  \cdot  
      \nonumber  \\
       &&\!\!\!\!{\small  \ba{ccccc} 
      \!(1\!-\!\mu_2\lambda_1) & & & & \\ & \!\!\!\!\!\epsilon_2  & & &\\ & & 
      \!\!\!\!\!\!(1\!-\!\mu_2\lambda_3) & &
      \\ & && \!\!\!\!\!\!(1\!-\!\mu_2\lambda_4) &   \\ & & & &\!\!\!\!\!\!\!\!\!\!\!1 
      \!\!\!\!\!\!\!\!\!\!\!\ea} 
      \cdot 
      \nonumber \\ 
      &&\!\!\!\!{\small \ba{ccccc} \!(1\!-\!\mu_3\lambda_1) & & & & \\ & 
      \!\!\!\!\!\!\!(1\!-\!\mu_3\lambda_2)  
      & 
      & 
      &\\ 
      & 
      &\!\!\!\!\!\epsilon_3 & &\\ & 
      && 
      \!\!\!\!\!(1\!-\!\mu_3\lambda_4)  &
      \\ & & & &\!\!\!\!\!\!\!\!\!\!\!\!\!\!1 \!\!\!\!\!\!\!\!\!\!\!\ea }  \cdot  \nonumber  \\   
      &&\!\!\!\!{\small 
      \ba{ccccc}   
      \!(1\!-\!\mu_4\lambda_1) & 
      & & & \\ & 
      \!\!\!\!\!\!\!(1\!-\!\mu_4\lambda_2)  
      & & 
      &\\ & 
      &\!\!\!\!\!(1\!-\!\mu_4\lambda_3) & &\\ & 
      && 
      \!\!\!\!\!\epsilon_4  &
      \\ & & & &\!\!\!\!\!\!\!\!\!\!1\!\!\!\!\!\!\!\!\!\ea}   \nonumber 
      \eq
  Note that in   general, the $m$-th diagonal entry 
   of the above product is given by 
   \be 
   \left[\prod_{k=1}^{N-1} \left(1-\mu_k \lambda_m 
   \right)\right]_{\!mm}\,=\,\;\; \epsilon_{m}\prod\limits_{\substack{k=1 \\ k\neq 
   m}}^{N-1}\left(1-  
   \frac{\lambda_m}{\lambda_k}  \label{prodNm2term}
   \right)
    \ee

  Now when $\A_i$ is the Laplacian matrix, except for $\lambda_N=0$, all its eigenvalues are positive. 
  Assuming 
  w.l.g. that    $\{\lambda_1<\lambda_2<\cdots <\lambda_{N-1}\}$, it is easy to  see that when $m>k$, 
   the  product  in (\ref{prodNm2term})
   will begin to involve increasingly larger scalars. As a result, in the best case scenario, when 
   $m=N\!-\!1$, 
   all these 
   numbers become larger than one in magnitude, and depending on the iteration, the magnitude of the 
    product 
   of  the  remaining   $N\!-\!2$ 
   terms in (\ref{prodNm2term}) can exceed that of  $\epsilon_{N-1}$.

\begin{figure}[htb]
                   	\epsfxsize 15cm \epsfclipon
                   	\begin{center}
                      \includegraphics[width=15cm,height=8.5cm]{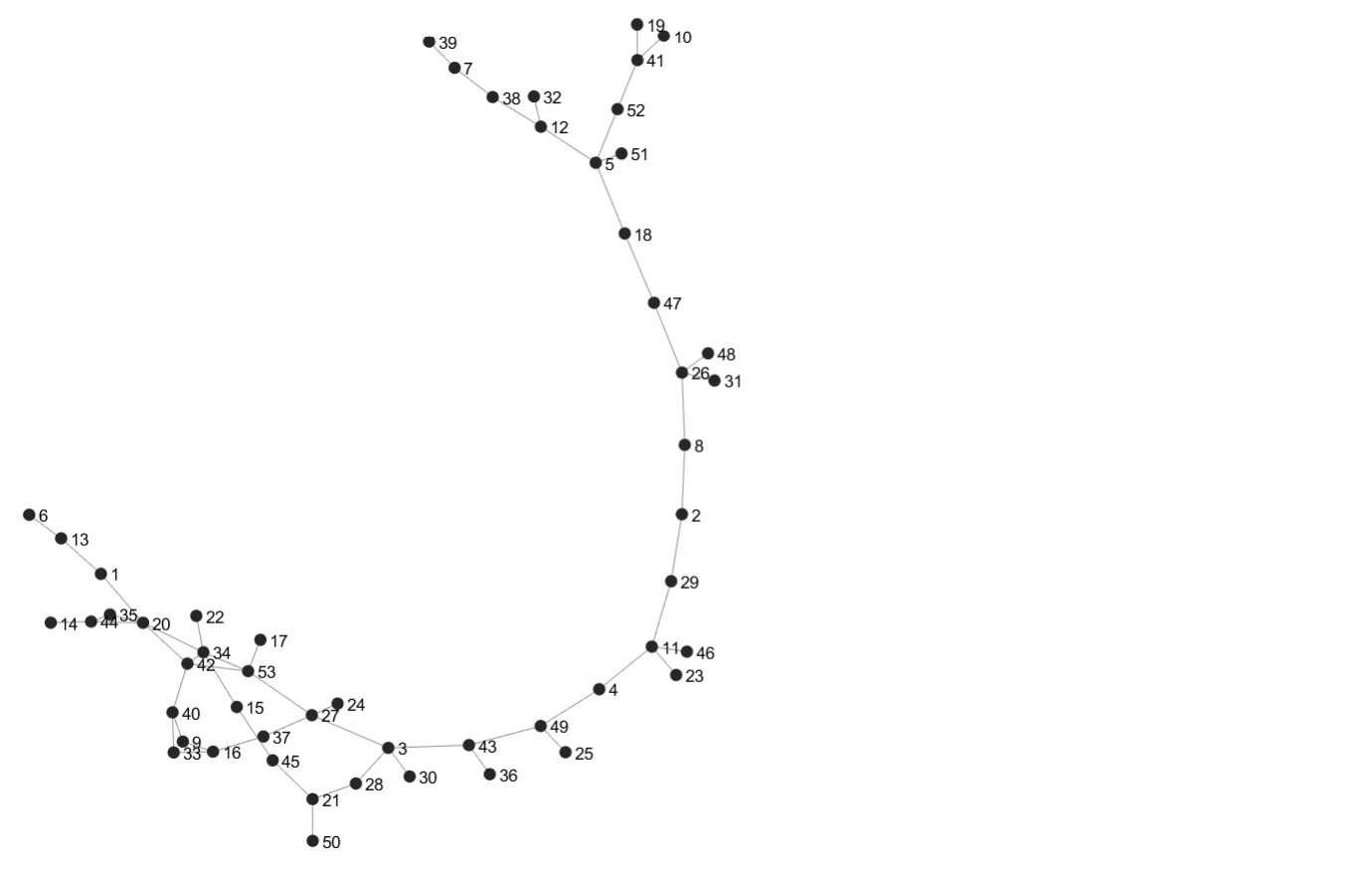}
                   		\vspace{-.5cm}
                   		  \caption{An undirected network with $N=53$ nodes, containing a linear  
                   		  strech sprouting from node 3 and exhibiting a few short 
                   		  branched paths.}\label{fig.netill}
                   	\end{center}
                   \end{figure}

   To illustrate this effect, we first consider the undirected  network of Fig.\;\ref{fig.net_asymm_20}, 
   for which we  compute the eigenvalues of $\C\tran\! \C={\LL}$,   taken as the 
      Laplacian 
      matrix (\ref{LLaplun}).    The range of nonzero eigenvalues of ${\LL}$ is $\{\lambda_{\sf 
    min}=0.021,\lambda_{\sf 
    max}=1\}$, where the smallest eigenvalue reflects the connectivity of the network. In this case, there is 
    no repeated eigenvalues, and it can be verified that the AC is achieved by the eigenstep method exactly  
    in $19$ 
    iterations. For 
    comparison, we also 
    considered the undirected graph     shown in 
   Fig.\;\ref{fig.netill},  
   comprising $N=53$ nodes, and where the respective nonzero eigenvalues of ${\LL}$  lie in the range 
   $\{\lambda_{\sf 
   min}=0.0088,\lambda_{\sf 
   max}=6.795\}$. 
                   
   Now, in   practice, when  performing  arithmetic operations,  some  error $\epsilon$  will naturally  
   arise  from 
   the   precision of the
   machine employed, i.e.,  it will take  place   with  ${\cal O}(\epsilon)$.
    For a specific hardware architecture, any floating-point computation will result in such error, say,     
    $\epsilon_c$,  which is
      defined as the smallest positive number such that $1+\epsilon_c \neq 1$. For example, Matlab's 
      precision 
   is such 
    that    
    $\epsilon_c=2^{-52} = 2.22\cdot 10^{-16}$\,(meaning the equivalent of  16 decimal 
   digits).  
   In  our context, this translates to an error for the $N$-th mode,
   \be 
   \epsilon_{N-1}\,=\,1-\mu_{N-1}\lambda_{N-1} \approx  2.22\cdot 10^{-16} 
   \ee
   It can be verified that  the corresponding product of the $N-2$ terms in (\ref{prodNm2term}) 
   results in $\approx 10^{42}$, 
   and moreover,    instability can be sensed as early as  $m=30$, when   (\ref{prodNm2term}) is 
   of the order of $ 
   10^{-2}$. This is already large, and  in this case the desired consensus will not be met.
    Note that in 
   this   example, the  effect of error accumulation  is the 
      same regardless of the ordering of $\{\lambda_k\}$. 
   
    \section{Discussion and Interplays} \label{discuss}
    
     We see that the complexity and difficulty in implementing 
    an 
   AC algorithm under the above constructions is first and foremost 
    dependent 
    on the constraint matrix $\C\tran \!\C$,  defined a priori.   For example, when $\A_i=I- \mu_i{\LL}$, 
    the number of 
    iterations required  to reach the exact  solution  via  
      eigensteps   depends on the 
      multiplicity 
      of the    Laplacian eigenvalues. Since  the 
            algorithm   
            requires  $K$ 
            iterations to converge\,(according to that 
            formulation),  in the worst case scenario it implies  a global  complexity of ${\cal O}(N^3)$ 
            operations\,(unless these  
                become untimely unstable   due to the 
            conditioning 
            of ${\LL}=\C\tran \!\C$).   Moreover, note from (\ref{QproQtra}) that the ordering of step 
            sizes $\{\mu_i\}$ for canceling  the modes of  $\A_i$ reveals to which of its  eigenvector 
            convergence will occur. For the Laplacian matrix, and according to (\ref{QproQtra}), the 
            surviving eigenvector will be $(\sfrac{1}{\sqrt{\!N}})\vones$. This suggests that the 
            convergence\,(and stability) may be affected considering a different ordering, since, e.g., 
            one    may  choose to cancel the  mode corresponding to the largest multiplicity  first. In the 
            same vein, we may order the  step sizes  in a way that another arbitrary eigenvector survives.  
            We   shall return to fill in   on these issues shortly, in the sequel.

      
     Note for now that, from the arguments that followed (\ref{brlsmnl0blms}), the eigenstep procedure holds 
     similarly for 
     arbitrary asymmetric matrices, not just Hermitian. Moreover,  it  would be   desirable that in this 
     case, 
     ${\LL}$ had 
     a high   multiplicity of eigenvalues, 
            so that 
            convergence is achieved faster.  On the other hand, if iterative algorithms are considered, one  
            would further  like 
            ${\LL}$ to exhibit the   smallest  eigenvalue spread possible.      These facts motivate not 
            only the use of optimized  matrices ${\LL}$,\,(and thus  
      $\A_i$), but   
      different  
      implementations, depending on whether the focus is on reduced complexity  or 
        improved   convergence.  \vspace{.15cm}

       \noindent {\bf Power Iterations.}    In a standalone scenario, the recursion 
        \be 
        \w_i = \A \w_{i-1} = (\I-\mu \A_o)  \w_{i-1} \label{shifpower}
        \ee 
 is simply  an (unnormalized) {\em shifted power iteration} \cite{Golub}. It is a widely known fact that   
 (\ref{shifpower})  converges  ideally
                              towards the scaled eigenvector   of   $\A$  corresponding to its 
                              larger 
                              eigenvalue in magnitude. This was in fact  predicted by  the exact 
                              solution\footnote{Recursion 
                                                             (\ref{genexpressfvdvv1}) relates to a  more 
                                                             general {\em subspace method}, 
                                                             in which the matrix iterations compute  several 
                                                             eigenvectors  at the same 
                                                             time; 
                                                              its normalized form  is referred to as  the    
                                                              {\em natural power 
                                                             iteration}\,\cite{Hua}.}  
                                                             (\ref{genexpressfvdvv1}). 
                               What is worth noticing is  that the convergence  of the power 
                               iterations to the  dominant eigenpair    holds similarly  when $\A$ is 
                               nonsymmetric, and not necessarily 
                               diagonalizable. 
                                  That is, assume 
                               for simplicity that $\A=\Q \D  \Q^{-1}$ is diagonalizable. Then, when $i 
                               \rightarrow 
                               \infty$,  we have 
                               \be 
                               \w_i\,=\, \lambda_N^i  (\u_1  \w_0)    \v_N \left(1 +  {\cal   
                               O}[(|\lambda_{N\!-\!1}/\lambda_N|^i)] \right)   \label{convOlamb}
                               \ee
                               where $\u_1 = {\sf e}\tran_1\Q^{-1}$, and which 
                               is in terms of the two larger  eigenvalues of $\A$ in magnitude, say, 
                                $\{\lambda_N,\lambda_{N\!-\!1}\}$.

                         Now,  consider the   roles played by the power method  in the 
                           pure  
                           linear algebra context, and in  the AC, distributed scenario. 
                           In the former case, the goal is to obtain the dominant eigenpair of a matrix, by 
                           {\em arbitrarily selecting} an initial condition $\w_0$, where   $\u_1\tran \w_0$ 
                           is 
                             unknown initially. 
                           Changing the value of $\mu$ will change the direction of 
                           convergence to a different eigenvector.  Note that $\w_i$ in (\ref{convOlamb})  
                           converges to the same form   of (\ref{genexpressfvdvv1}), except for a scaling 
                           factor, which grows or 
                              decreases with time, unless $\lambda_N=1$. 
                            This is in fact one source of 
                              instability  observed in the eigenstep method: no normalization is employed, 
                              and as a result, any numerical error is bound to accumulate.  In the standard 
                              power iterations, this is resolved 
                              via 
                         normalization of $\|\w_i\|$\,(either by  the   Euclidean or infinity  norm, which 
                         can be duly
                         justified in 
                         an optimal sense), so that 
                         $\w^o \longrightarrow   \alpha \v_N$, with    $\alpha = \pm 1$.   The important 
                         fact 
                         to note    is that  convergence to the dominant eigenpair occurs {\em 
                         regardless} of the   structure of $\A$\;[except either when 
                             $\{\lambda_N,\lambda_{N\!-\!1}\}$  have the  same magnitude,   if they are too 
                             close, or if  $\w_0$ is     chosen (unlikely)  orthogonal to  $\u_1$, as can be 
                             verified from 
                             (\ref{convOlamb})].
                             
                             In the distributed context, 
                             on the other hand, 
                           $\w_0$ 
                           is {\em given 
                           but unknown}. An asymmetric matrix  $\A$ in this case  may correspond 
                          either  to a  directed or an   undirected 
                           graph\,[see, e.g., (\ref{normL})],  which {\em is 
                           entirely    known, including its 
                           eigenstructure}, and which may comprise    
                            negative or  even complex eigenvalues/vectors. 
                                       From (\ref{convOlamb}), we see that the AC can be obtained by 
             replacing $\w_0$ by    $\w_0'= \D_1^{-1}\w_0$, and by computing
                   \be 
                    \w^o\,=\, \D_2^{-1}\w_i\,=\, 
                      \dfrac{1}{N}\vones \!\vones\tran \w_0
                     \ee
                      where 
                     \begin{subequations}
                      \bq 
                      \D_1&\!\!=&\!\!(\sfrac{1\!}{\sqrt{\!N}})\text{diag}(\u_1)  \label{D2diaff0}  \\
                      \D_2&\!\!=&\!\!(\sfrac{1\!}{\sqrt{\!N}})\text{diag}(\v_N) \label{D2diaff}
                      \eq
                      \end{subequations}
                  Observe that in this case, we need to assume that $\{\v_{\!N},\u_1\}$ have nonzero entries, 
                  and that their   dynamic range is properly accounted for. That is, many 
                     entries of $\{\v_{\!N},\u_1\}$ may assume very small numbers, specially when $\A_o$ is 
                     designed to have 
                     small $|\lambda_{N\!-\!1}|$\,[The usual approach is  to constrain  
                     $\A$ to have its maximum eigenvalue associated to $\vones$, as stated in  
                     (\ref{probAfor})]. 
                    Still, any (sparse) matrix possessing $\{\v_{\!N},\u_1\}$ with    nonzero  
                                entries,  
                                and 
                              with $|\lambda_N|$ sufficiently larger than $|\lambda_{N\!-\!1}|$, can 
                              serve as a solution to the AC   problem.

                           The important aspect here is to realize  that,   even though normalization   
                           implies 
                           computation 
                           of a global information,  $\A$  is known beforehand; it can be  
                           designed to 
                           have $\lambda_N=1$, and  moreover, possibly having  the ratio 
                           $|\lambda_{N\!-\!1}|/|\lambda_N|$ 
                           minimized  for optimized  convergence rate. For example, 
                           starting from a {\em normalized} Laplacian matrix,  we can restrict its spectrum 
                           to 
                           $0=\lambda_1 < \lambda_2 <\cdots   < \lambda_{N} < 
                             2$ by defining
                             \be 
                              \xoverline{\LL}\,=\, ({\cal D}^{\dagger})^{\!1/2} {\LL} 
                              ({\cal D}^{\dagger})^{\!1/2}\;\;\;\;\text{or}\;\;\;\; 
                              \xoverline{\LL}\,=\, {\cal D}^{\dagger} {\LL}    \label{normL}
                             \ee
                          where ${\cal D}=\text{diag}({\LL})$ is the degree matrix\,\cite{Chung}. As a 
                          result,  the matrix $\A_o = \I -  
                                                     \xoverline{\LL}$ will have $\lambda_N=1$\,(Note 
                          that the 
                          right-most 
                          expression above implies that the combination matrix of 
                          its undirected network is not symmetric).

 The above discussions suggest  that   the  eigenstep procedure is not   limited to undirected networks,   
 but can be applied similarly  to any directed graph, 
 given    proper 
 normalization.      That is, consider any 
 diagonalizable matrix  $\A_o=\Q 
 \bLambda 
 \Q^{-1}$, where $\bLambda$ has the form (\ref{lmultipl}).  Instead of  the shifted power iteration 
 (\ref{shifpower}), consider the  recursion with time-varying coefficients
            \be 
                \w_i = \left(\I - \dfrac{1}{\lambda_i}\A_o\right) \w_{i-1} \;\;\;i=1,2,\ldots,K 
                \label{shifpowertv}
              \ee   
       Solving for $\w_i$, we have
    \bq 
     \w_K &=& \Q \prod_{1=1}^{K}\left(I - \dfrac{1}{\lambda_i}\bLambda \right) \Q^{-1} \w_{0} \nonumber  \\
     &=& \Q \ba{cccc} 0 &&& \\ &  0 &&\\ && \ddots & \\ & && \!\!\prod\limits_{i=1}^{K}\!\left(1-  
        \frac{\lambda_N}{\lambda_i}\right)  \ea \Q^{-1} \w_{0}  \nonumber \\
     &=& \prod\limits_{i=1}^{K}\!\left(1-  \frac{\lambda_N}{\lambda_i}\right)  (\u_N \w_0) \v_N 
     \label{scaledeig}
    \eq       
    where $\u_N$ is defined  in (\ref{uNQinv}).     Hence, similarly to  the power iterations, the 
    convergence is towards the scaled eigenvector 
    $\v_N$.   The scaling factor in (\ref{scaledeig}) can then be 
     removed 
     by normalizing  
    (\ref{shifpowertv}) by $(1-\lambda_N/\lambda_i)$ so that 
    \be 
    \A_i =  \dfrac{1}{\lambda_i \!- \!\lambda_N} \left(\lambda_i \I -  \A_o\right) \;\;\;i=1,2,\ldots,K 
    \label{AiK111N}
    \ee 
    
    Note that when $\A_o={\LL}$, i.e., the 
             Laplacian 
             matrix, its eigenvalues can be ordered so that $\lambda_N=0$, and  (\ref{AiK111N}) collapses to 
             (\ref{weightsAil}).     
             Note also that in theory, the choice of $\A_o$  is  not limited to the Laplacians; any matrix, 
             which can have complex,  or even a minimized number of elements\,(targeting minimized 
             complexity), can be considered. Since the eigenvalues $\{\lambda_i\}$ can be 
             negative, there is more room for numerical errors to be reduced, given that the denominator in 
             (\ref{AiK111N}) may 
             assume larger values, compared, e.g., to the caculation in (\ref{prodNm2term}).  In addition, 
             notice that  
                the error 
                          from the AC may be reduced at a faster rate   than $K$ iterations  for matrices 
                          with 
                          higher 
                          multiplicity of eigenvalues, at the expense of 
                          increased dynamic range in data representation. 
        That is, defining 
          \be
           \D_1\,=\,(\sfrac{1\!}{\sqrt{\!N}})\text{diag}(\u_N) 
              \ee                  
            and $\D_2$ as in  (\ref{D2diaff}), with  $\w_0'= \D_1^{-1}\w_0$, the AC  becomes
                                          \be 
                                             \w^o\,=\, \D_2^{-1}\w_i\,=\, \dfrac{1}{N}\vones \!\vones\tran 
                                             \w_0
                                           \ee  
        where $\{\v_N,\u_N\}$ may have small entries. 
  Observe that all we have said regarding error accumulation in (\ref{prodNm2term})  still holds,  its effect 
            depending  on $\A_o$  and its size. Yet, a different choice for $\A_o$ may be sufficient to 
            counter the 
            numerical 
            errors that arise, depending on the magnitude of  $\{\v_N,\u_N\}$.

               \noindent {\bf Accelerated Gradient Algorithms.} Most AC strategies in the 
               literature 
               assume 
               that  the
              network iterations  are of the linear form  $\w_i= \A_i \w_{i-1}$.   
              %
               This approach is not unique, and in fact, given a Laplacian matrix ${\LL}$, 
                there is a 
                variety of 
                accelerated    algorithms that can be 
                      considered, including the 
                      Nesterov accelerated 
                      gradient\,(NAG)\,\cite{Nesterov},  
                      Adaptive 
                      Subgradient\,\cite{Duchi}, and Chebyshev semi-iterative methods\,\cite{Golub1}. 
                        As we have  discussed, because normalization is generally a sensitive step, an 
                        alternative 
                        approach is to 
                           reconsider 
                           fixed step sizes  under these alternative recursions. This avoids normalization by 
                            global parameters that require online computations, so that the node estimates 
                            can 
                            be  
                            combined and updated  independently.  For the sake of comparison,   here we 
                            consider 
                        one 
                      possible variation of the NAG recursions\,\cite{Nakerstn}, 
                       applied to  the newtork context:
                      \begin{subequations}
                      \bq 
                        \q_i &=&  \beta \q_{i-1} +  \alpha  {\LL}(\w_{i-1}-\sigma \q_{i-1})  
                        \label{NAGa} \\
                          \w_i &=& \w_{i-1} - \q_i \label{NAGb}
                     \eq
                    \end{subequations}
                     where $\{\alpha,\beta,\sigma\}$ are pre-defined learning parameters.
                     

                       Although the structure of the above  recursions   taps on one 
                     of the oldest problems of optimization, the advantage here is  that these parameters are 
                     combined and 
                     adapted {\em 
                     entrywise}, without  covariance inversion.        
                      Moreover,  since ${\LL}$ is fixed, there is no need for the nodes to 
                     store an arbitrary number of coefficients, given that  the time-varying  nature of the 
                     above recursions is   only 
                     due 
                     to their coupling.

\begin{figure}[htb]
                                             	\begin{center}
                                             		%
                                             	 \includegraphics[width=9.5cm]{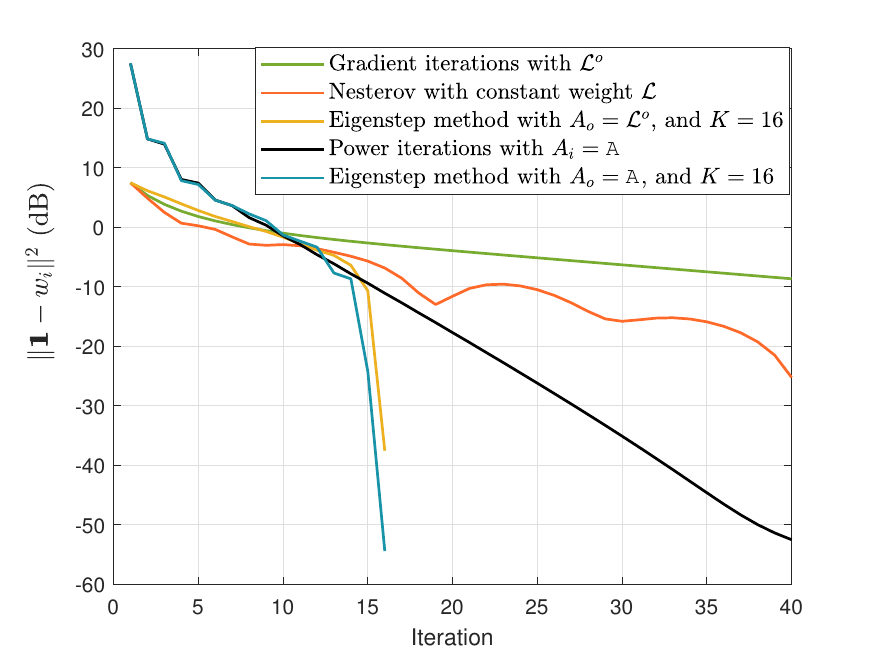}
                                             		\vspace{-.5cm}
                                             		  \caption{Convergence of  several methods  to the 
                                             		  exact  solution  for the undirected network of 
                                             		  Fig.\;\ref{fig.net_asymm_20}\,(with target vector 
                                             		  normalized to $\vones$).}\label{fig.eig_x_power_20sym}
                                             	\end{center}
                                             \end{figure}

                      At this point, it is worth  confronting the performance of 
                     the 
            eigensteps  with the ones of the power recursions that make use of arbitrary and optimized 
            matrices, as 
            well as 
            with (\ref{NAGa}-b) in the distributed case.   We consider initially  the convergence  
            over    the undirected graph of Fig.\;\ref{fig.net_asymm_20}, shown in   
        Fig.\;\ref{fig.eig_x_power_20sym}, where the distances   between  $\w_i$ and   the (normalized) 
        target vector, $\vones$, are  
                  computed.       First, we  observe that the worst performance is due to the gradient 
                  iterations  
                  in terms 
                  of the 
                          optimized Laplacian
                          matrix, i.e., 
                  $\A_i=\A = \I- {\LL}^o$,  which is obtained by solving  (\ref{probAfor}).  
        We denote the eigenvalues of ${\LL}^o$, as well as the ones of the adjacency matrix ${\sf A}$, as 
        listed below:\vspace{.1cm}
        \be 
 {\small  \begin{array}{c} \hspace{-2.5cm} \text{\normalsize eig}({\LL}^o) = \\   \{ 
            \stackrel{\phantom{c}}{0.000000000000000}, \\ 
          \;\,0.052535234682282,\\
           \;\,0.052535264553794,\\
           \;\,0.083436170825888,\\
          \;\, 0.181855637792896,\\
          \;\,0.401856564107921,\\
           \;\,0.427189140429399,\\
          \;\, 0.495504631874321,\\
          \;\, 0.686002656764135,\\
          \;\, 0.925663817734459,\\
          \;\, 0.999999924855551,\\
           \;\,1.000000005771190,\\
          \;\, 1.000000234186339,\\
           \;\,1.171675928308657,\\
          \;\, 1.183654232753852,\\
           \;\,1.259104383988848,\\
          \;\, 1.397001250816675,\\
          \;\, 1.762431475613161,\\
          \;\, 1.947464654320199,\\
         \;\,1.947464756759985\}  \end{array}}  \;\;\; \;\;\; 
           {\small  \begin{array}{c} \hspace{-2.5cm} \text{\normalsize eig}({\sf A}) =\\  \{ 
                       \stackrel{\phantom{c}}{-2.500705872939440}, \\
                      \;\, -2.209324772072598, \\
                       \;\, -1.855625375629181, \\
                      \;\,  -1.530096384802822, \\
                      \;\,  -1.396508043500946, \\
                      \;\,  -1.146042533665765, \\
                      \;\,  -1.000000000000000, \\
                      \;\,  -1.000000000000000, \\
                      \;\,  -1.000000000000000, \\
                       \;\, -0.674484311935240, \\
                      \;\,  -0.550540360935937, \\
                       \;\, -0.356647333315697, \\
                        \;\, \;\;\;0.000000000000000, \\
                       \;\,  \;\;\;0.125780164164776, \\
                       \;\,  \;\;\;0.650304387934893, \\
                        \;\, \;\;\;1.253890737507349, \\
                        \;\, \;\;\;2.277240308122007, \\
                        \;\, \;\;\;3.078461594541966, \\
                       \;\,  \;\;\;3.383236084445506, \\
                        \;\, \;\;\;4.451061712081128\} \end{array} }  \nonumber 
         \ee
         
         We see that convergence in this case is a leisurely process, and is  followed by   the Nesterov 
         recursion 
         (\ref{NAGa}-b), which   was tunned offline  with $\{\alpha=0.15,\,\beta=\sigma=0.85\}$ for 
 fastest convergence. 
         The 
         standard power iteration 
         on the 
         other hand, is 
         implemented with the adjacency matrix
         $\A_i={\sf A}$,  which further improves convergence over the latter. This comes 
         at 
         the expense of an increase in dynamic range by  $\approx 100$, due to the small entries of the 
          dominant 
         eigenvector $\v_N$. The eigenstep was then tested, first  with  an optimized Laplacian, $\LL^o$,  
         exhibiting 2 eigenvalues with 
         multiplicity 2, and 1 eigenvalue with multiplicity 3\,(i.e., roughly, since these are 
         numerically very close), and the rest with multiplicity 1. Although in theory  
         convergence should 
         occur 
         in $K=15$ 
         iterations, since the double multiplicity of the largest eigenvalue of ${\LL}^o$ is not 
         exact, we assumed  that $K=16$. The eigenstep method implemented via (\ref{weightsAil}) has  
         $\v_N=(\sfrac{1\!}{\sqrt{\!N}})\vones$ and 
         attains approximately 
         $-40$\,dB at $i=16$, which, although stable,  still exhibits some 
         residual error. For 
         comparison, we considered the eigenstep method that makes use of  $\A_o={\sf 
         A}$, via the  normalized recursion (\ref{AiK111N}). The eigenvalue 1 has multiplicity 3 in this 
         case, and even under a larger dynamic range\,($\approx 100$), it leads to a smaller error
         (in addition to
         being   less complex, since most entries are unity). 

            The performance of   the  eigenstep method for a   directed graph\,(with connecting structure   
            listed   in the Appendix A),  is depicted in Fig.\,\ref{fig.eig_x_power_20}. In this example, we 
            set the combination matrix as  
          a (right) Laplacian,    
          in   that ${\LL}\vones=0$, however, asymmetric. The ordering of the stepsizes is such 
          that 
          the surviving 
          eigenvector is  $\v_N=(\sfrac{1\!}{\sqrt{\!N}})\vones$\,(red curve). In this case, convergence 
          is attained exactly in $N-1=19$ iterations,  
          which, as before, 
          is not 
          gradual, but undergoes a hasty drop   at $i\!=\!19$. For comparison, we designed a combination 
          matrix 
          $\A_o=\A_c$, 
          with 
          $\lambda_N=1$, and whose second largest eigenvalue has magnitude $\lambda_{N\!-\!1}=0.37$, in order 
          to keep the dynamic range  of  $v_N$  $\approx 100$, as before.  It is plain to see
          that convergence of the power method in this case is  faster, and although   not exact, 
          it reaches $-50$\,dB in $19$ iterations\,(black curve). Using this same matrix for the 
          eigen-iterations, convergence is not only faster, but also exact\,(dark yellowish curve). 
          
           \vspace{-.2cm} 
              \begin{figure}[htb]
                                   	\epsfxsize 9.65cm \epsfclipon
                                   	\begin{center}
                                   		 \includegraphics[width=9.5cm]{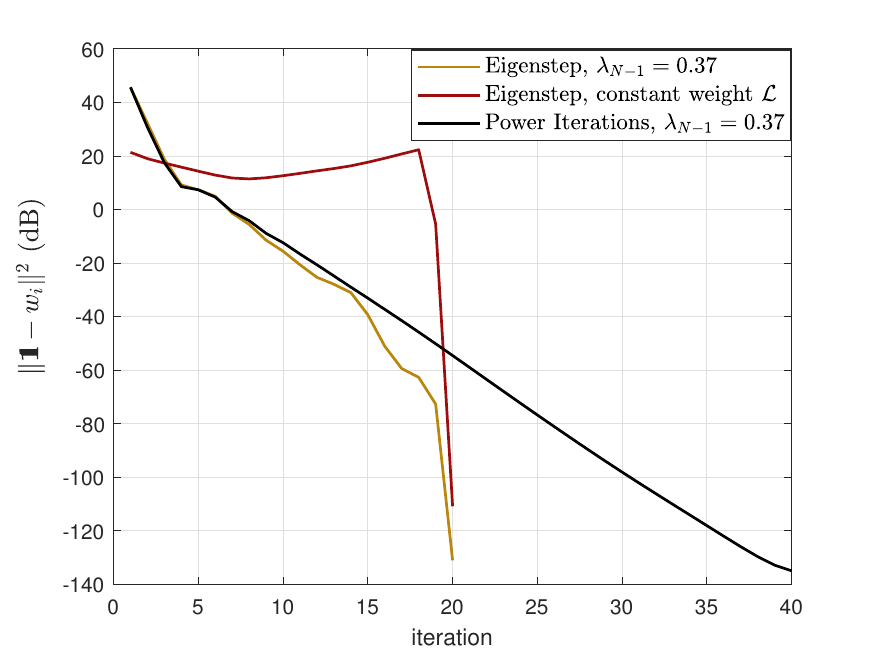}
                                   		\vspace{-.5cm}
                                   		  \caption{Convergence of the eigenstep method  for
                                   		  ($i$)\,$\A_o={\LL}$, ($ii$)\,$\A_o=\A_c$, and ($iii$)\,power 
                                   		  iteration with 
 $\A_i=\A_c$,   for the  directed 
graph of   Appendix A. In the former case, the curve  lies flat until sudden, exact 
convergence.}\label{fig.eig_x_power_20}
                                   	\end{center}
                                   \end{figure}

\begin{figure}[htb]	\epsfxsize 14cm \epsfclipon 
                             \begin{center}
                       \phantom{xxx}\includegraphics[width=0.85\textwidth]{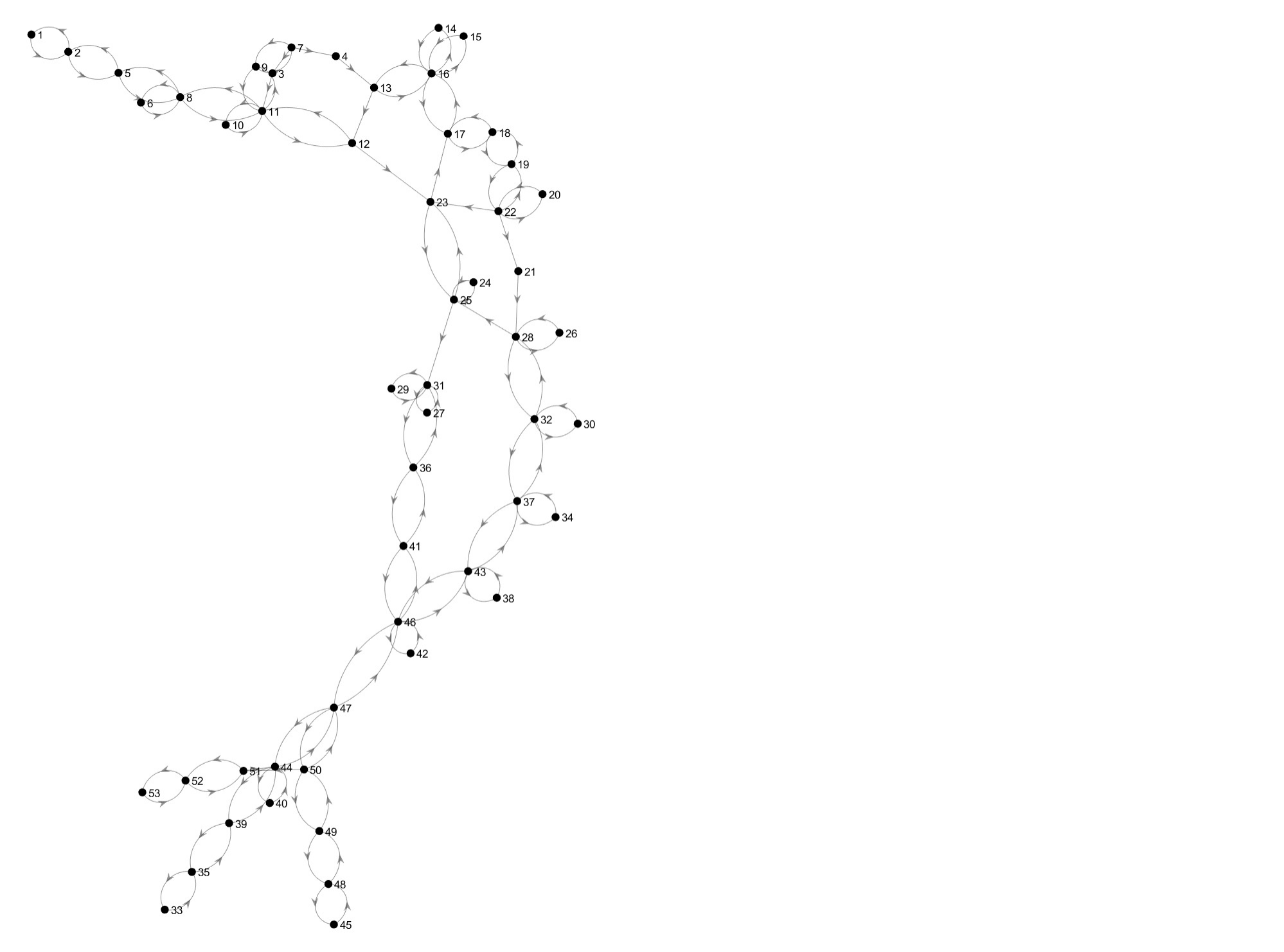}
                      	\vspace{-.6cm}
                             \caption{Directed graph constructed   from  the  network of   
                             Fig.\;\ref{fig.netill}.}\label{fig.netillu} 
                             	\end{center}
                           \end{figure}

  The behavior of these algorithms for  a larger graph  is not necessarilly the 
  same. We  illustrate  this fact by reconnecting  some of the edges of the 
           network of 
              Fig.\;\ref{fig.netill} to farther  nodes, in order to construct a directed topology,    
              depicted 
              in 
              Fig.\;\ref{fig.netillu}.  \\
               
\begin{figure}[t]
                                             	\epsfxsize 9.4cm \epsfclipon
                                             	\begin{center}
                                             		\includegraphics[width=9.5cm]{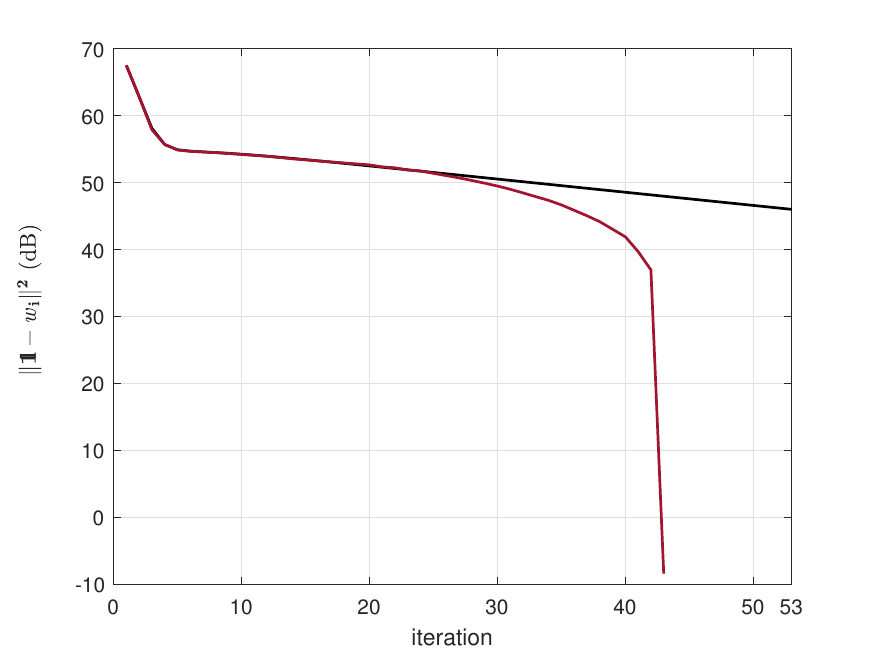}
                                             		\vspace{-.7cm}
                                             		  \caption{Convergence of the eigenstep method to the 
                                             		  exact 
                                             		  solution in 
                                             		  $N=43$  
                                             		  iterations, when $A_o={\sf A}$, i.e., the adjacency 
                                             		  matrix.}\label{fig.conv_eig}
                                             	\end{center}
                                             \end{figure}

              Here, the 
              resulting matrix $\A_o$ is 
              asymmetric, which we set, for this example, as the  adjacency matrix,   
              $\A_o={\sf A}$. 
              Again, the overall 
                       computation 
                       requires mostly   additions, except for the scaling factors in (\ref{AiK111N}) 
                       applied 
                       to each node 
                       individually.  
                   Figure \ref{fig.conv_eig} illustrates the convergence to the vector $\vones$ 
                  after $N=43$ 
         iterations.  Note that by virtue of (\ref{Alkcond}), the 
                     minimum number of iterations $j_o$   for the recursions to achieve AC is $j_o=20$. 

We sum up the main aspects of this discussion as follows:

      \vspace{-.2cm}
   \begin{enumerate}
   \item  In theory, consensus can be achieved via the eigensteps or standard power iterations, for 
   undirected and 
   directed 
   graphs, where the adjacency matrix is asymmetric;  
    \vspace{.07cm}
   
   \item  In principle, the structure of $\A_o$ is arbitrary, and its coefficients can be generally complex. 
   This 
   means that   AC is attained by   
   pre- and post-multiplication  at the nodes   according to (\ref{D2diaff0}-b);   \vspace{.07cm}
   
    \item  The eigenstep method becomes faster depending on the multiplicity of the eigenvalues of the 
    underlying combination matrix.  The Laplacian matrix will not necessarily yield the highest multiplicity  
  compared, e.g., to the one of a directed graph.
 In 
    Sec.\,\ref{sec.exact},  we shall  minimize the number of elements of the Laplacian matrix for this same 
    graph, 
    leading to a multiplicity 6 for the eigenvalue at $1$;    \vspace{.07cm}
   
   \item  The convergence of the eigenstep iterations is not dependent on the 
   magnitude of the second eigenvalue, which could be the same as the largest one.  Since all 
   eigenvalues are known in advance, in theory convergence is achieved exactly in at most $N\!-\!1$ 
            iterations. Also, since these  are not unique,   alternative matrices may yield 
            different numerical behaviors and faster error reduction.   Still, 
            error accumulation 
            may or may not be larger in practice than what would be observed via the   power method, as seen 
            in the above examples.  One can in 
   fact  design combination matrices such 
       that   $|\lambda_{N\!-\!1}|$ is minimized, improving 
         convergence of the power iterations, which become much faster   and more exact than the eigenstep 
         method, at the expense of increased dynamic range of the nodes information. 
 \end{enumerate}

   \section{Exact Computation of the solution} \label{sec.exact}
   
  The numerical difficulties that arise  in all the above constructs can be tackled    by 
  envisaging an exact 
  implementation 
  that does not rely on  gradient-like iterations\,(which in the case of 
  ill-connected networks  might take a long time to converge). We shall address this issue by 
  considering first an undirected network.  
  Here, any constraint matrix  $\C$  satisfying $\C\w=\0$ will do, and 
        the 
        goal is to 
        implement   (\ref{brlsmnl0002})  as efficiently as possible, 
        with   focus on the design of $\C$ with an eye on complexity reduction.   We  
        propose two 
        solutions 
        that  may take either the same number of iterations of the 
        eigenstep method, or at most $2N$ 
        iterations, in a sequential manner. The overall complexity will be considerably reduced to ${\cal 
        O}(2N)$ 
        {\em additions} in the latter, in contrast to ${\cal O}(a N^3)$ {\em multiplications} and additions 
        for 
        the eigenstep,  
        or some iterative solution. 
   
   Thus, let us return to the exact solution    (\ref{brlsmnl0002}) (with $\f=\0$):
   \be
  \w^o \,=\, \Big[\I -    \C\tran \!\left( \C   \C\tran\right)^{\!-1}  \!\C\Big]  \w_0   \label{exctsolwo0}
   \ee
   To motivate the idea,   recall that one  immediate approach to solving a  linear system of 
   equations  
   $\LL\w= \C\tran \!\C\w =\0$, is to express it  in terms of  triangular factors, say, 
   $\L$\,($\L\tran$), 
   so 
   that operations with these matrices 
   can be solved 
   efficiently via successive forward and back-substitutions\,(this is similar  to a 
   Cholesky decomposition  of positive definite Hermitian matrices, although here ${\LL}$ is only positive 
   semi-definite).   Still, 
   we can 
     manage to decompose $ {\LL}$ as  ${\LL}=\L \L\tran$, where $\L$ is no longer unique, but 
   generally a tall   matrix.   More specifically,  we would like to  construct a   matrix 
   ${\LL}$ such that 
   \be 
    {\LL}\w\,=\, \L \L\tran \w\,=\,\0 
    \ee
    where $\L$ is a $N \times (N-1)$ sparse, lower triangular, full column rank matrix. Since $\L \L\tran$ 
    and 
    $\L$ 
    have the same nullspace, we simply require that  
    \be 
    \L\tran \w=\0  \label{LTeq}
    \ee
    
    In order to conform  with the 
     network structure, (\ref{LTeq}) can be  satisfied by $(i)$ designing  the columns of  $\L$ such that 
     they 
     add up to 
     zero, $(ii)$   having 
   its lower 
   triangular part exhibiting a zero pattern that coincides with the one of the adjacency matrix, and 
   $(iii)$\,selecting the main diagonal elements of $\L$ to be positive. The latter is necessary to ensure  
   that, by 
   the 
   Sylvester's law of inertia, ${\LL}$ will have $N-1$ positive eigenvalues and a simple eigenvalue at 
   zero. 
   
\begin{figure}[htbp]
   	\epsfxsize 9cm \epsfclipon
   	\begin{center}
 \epsffile{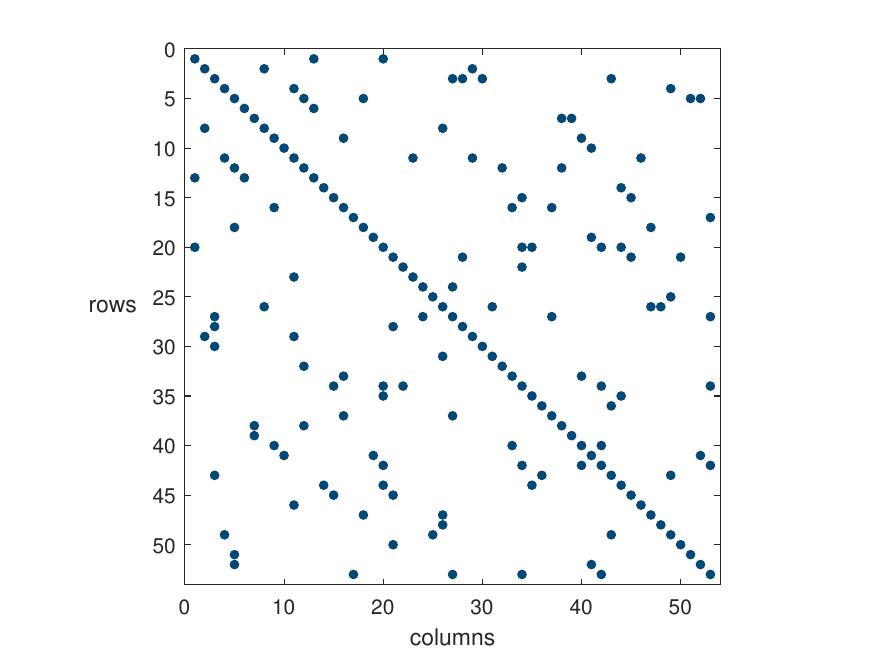}
   		\vspace{-.5cm}
   		 \caption{Support of the adjacency matrix\,(with self-loops) for the $N$-node graph of 
   		Fig.\;\ref{fig.netill}.}\label{fig.Adjmat}
   	\end{center}
   \end{figure}

   A natural question  is whether it is possible to construct $\L$ in this way, starting from an 
   arbitrary 
   graph 
   with a predefined node  numbering. To elaborate on this 
   point, we continue with the example 
   of Fig.\;\ref{fig.netill}, for which its support\,(including self-loops), denoted by  ${\sf A}$,  is 
   depicted 
   in Fig.\;\ref{fig.Adjmat}.      
   Assume that  we attempt to construct an $N \times (N-1)$ (Laplacian-like) lower triangular factor $\L$ 
   from 
${\sf 
   A}$, i.e., by defining its support   as
   \be 
   \A_L \,=\, \left[\, \text{lower}\{{\sf A}\} \,\right]_{1:N,1:N-1} \label{ALlow}
   \ee
   In this case, several columns will end up having a single unity entry, so that their sum cannot be made  
   zero. 
     Similar issues, which commonly arise  in related  sparsity problems, are dealt with by first renumbering 
     the 
    nodes of the network according to some pivoting scheme. The idea is to cluster the nonzero entries of 
    each 
    row  
    around 
    the main diagonal so that its profile can be  reduced. This is in turn equivalent to applying a 
    permutation 
    matrix 
    $\PP$ 
     to ${\sf A}$, while reordering the nodes accordingly. The procedure then yields a new adjacency 
    structure, 
    say, $ {\sf 
    A}_p$:
    \be 
    {\sf A}_p\,=\,  \PP {\sf A}  \PP\tran  \label{reordPm}
    \ee
   For this purpose, several techniques exist which tend to narrow down the ``bandwidth'' or the ``envelope'' 
   around the main diagonal, 
   including   
   the {\em Symmetric Approximate Minimum Degree Permutation}\,(SAMD)\,\cite{SAMD}, the {\em Nested 
   Dissected}\,(ND) 
   method\,\cite{Karypis}, 
   and the
   {\em Symmetric Reverse Cuthill-McKee Permutation}\,(SRCM) algorithm\,\cite{George}. We shall focus on the 
   latter, which 
   aims 
   at   
   minimizing the distance between a   vertex and its ordered neighbors.  This step  is 
   usually performed    prior to a Cholesky or LU decomposition in solving linear equations, and can be 
   found, 
   e.g.,  as 
   a   MATLAB function.
   
\begin{figure}[htbp]
    	\epsfxsize 9cm \epsfclipon
    	\begin{center}
   	\hspace{-0.3cm}\epsffile{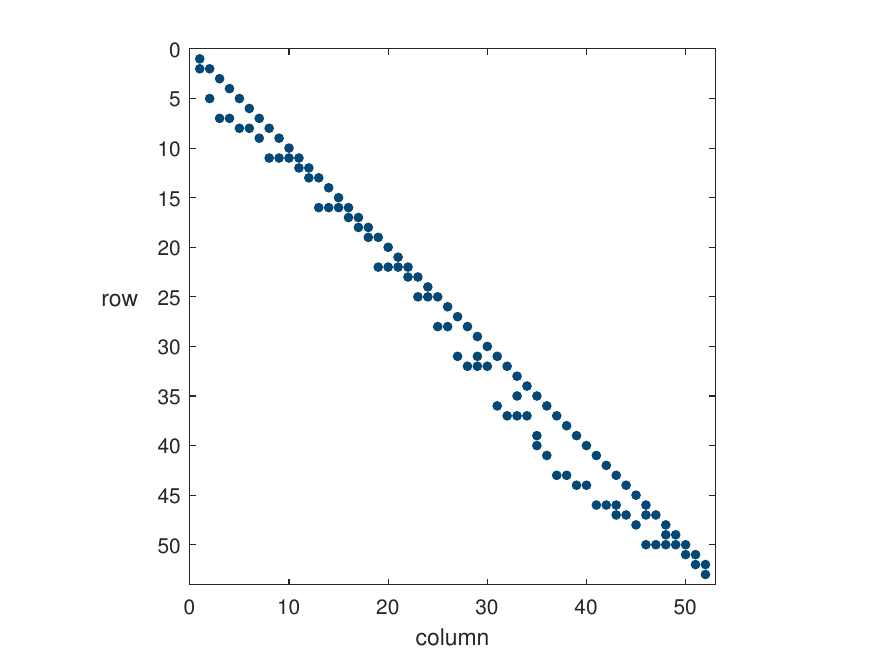}
    		\vspace{-.1cm}
     \caption{Lower triangular support of $\L$ after reordering via the SRCM 
    		method.}\label{fig.Adjmat1}
    	\end{center}
    \end{figure}

   Now, it is a  known fact that   after the application of the SRCM permutation algorithm,   the reordered 
   structure of 
   ${\sf 
   A}_p$  in a connected network 
   becomes  block 
   tridiagonal\,(see, e.g., \cite{Pedroche},\cite{Reid},\cite{Saad}).    
   In this case, the corresponding lower triangular factor $\A_L$ in (\ref{ALlow}) will end up having the 
   desired lower triangular structure\,(or any 
   procedure that succeeds in structuring the adjacency profile in this way would do).  For the  
   graph example above, reordering  results in the profile for 
   \be 
   \A_L 
   \,=\, \left[\, 
   \text{lower}\{{\sf A}_p\} \,\right]_{1:N,1:N-1} \label{reordprof}
   \ee 
     illustrated in Fig.\;\ref{fig.Adjmat1}. Observe that  
     each column of $\A_L$ has at least one entry below the one  in its main diagonal.  
    
   A simple way to construct the entries of $\L$ is to define   
    \be
    [\L]_{k\ell} \;=\;\left\{ \begin{array}{ll} \!\!1, & \text{if }k=\ell \\\!\!-1/(n_k\!-\!1),  & \text{if 
    }k\!>\!\ell  
    \text{ with nodes $k$ and $\ell$ $\in {\cal N}_k$} \\ \!\!0, & \text{otherwise} \end{array}   \right. 
    \label{LLapl}
    \ee
    where the  column-sum  of $\L$ now conforms with 
     the   one   of the normalized 
      Laplacian matrix defined in (\ref{normL}). 
Note that here  we do not require ${\LL}$ to have the same profile of  ${\sf A}_p$. The former can 
in 
    fact  exhibit fill-in entries compared to the latter, and what is actually required  is that the 
   {\em solution}
    conforms 
    with the network topology\,(not the statement of the problem).  
     
 As a result, following (\ref{exctsolwo0}),  expression (\ref{LTeq}) is now in terms of $\L$, and given by 
    \be 
    \w^o \,=\,  \Big[\I -    \L \!\left( \L\!\tran \!\L\right)^{\!-1}  \!\L\tran\Big]  
    \w_0   
    \label{exctsolLL}
    \ee
    where $\L$  can be partitioned as
    \be 
    \L\,=\,\ba{c} \bar{\L} \\   \b \ea 
    \ee
    This expression has the same form of (\ref{DecQDQsgen}), which leads to (\ref{prodmatswog}) with 
    $\C_1=\C_2=\L^{\sf T}$:
    \bq
    \w^o &\!\!=& \!\!\gamma \underbrace{\ba{cc} \!\bar{\L}^{\sf -T} & \0 \\ \0 & 1 \ea}_{ \W^{\;\sf  -T}} 
    \underbrace{\ba{c} 
    \!\b\tran \\ 
    -1\ea}_{\b_1\tran}   \!\ba{cc} \!\b & -1 \ea  
    \!\ba{cc} \!\! 
    \bar{\L}^{\sf -1} & \0 \\ \0 &  1 \ea  \w_0  \nonumber   \\
   &\!\!=&  \!\!\gamma  \W^{\;\sf  -T}  \b_1\tran \b_1 \W^{\;-1}   \w_0\label{prodmatswo}
    \eq


    From this expression, we can promptly devise two alternative  distributed AC implementations, as we now 
    explain. 
    \vspace{0.16cm}

    \noindent  {\sfb 1)  \underline{Network back-substitution implementation}}. \vspace{.1cm}
     
    Consider the   vector 
    partitions
    \be 
      \w_0 = \ba{c}  \w_0(1\!:\!N\!-\!1) \\ w_0(N) \ea,\;\;\;\;   \w^o = \ba{c}  \w^o(1\!:\!N\!-\!1) \\ 
      w^o(N) \ea
    \ee
    Starting with $\w_0$,  we first  solve the triangular system of equations via back-substitution:
    \be 
    \bar{\L} \w_1(1\!:\!N\!-\!1) \,=\,   \w_0(1\!:\!N\!-\!1) 
    \ee
    Since $\bar{\L}$ conforms with ${\sf A}_{p}$, the operations are distributed, and the same applies to 
    the 
    multiplication with the vector $\b$, defined through (\ref{prodmatswo}). Now, define  the scalar quantity
    \be 
    w_1(N)\,=\,\gamma  \left[ \b \w_1(1\!:\!N\!-\!1)- \w_0(N)\right]\label{wiw0w0w9}
    \ee
    which will be available at node $N$ at time $N$. Then,
    \be 
    \w^o\,=\,  \ba{c} \bar{\L}^{-{\sf T}}\b\tran w_1(N)   \\ -w_1(N) \ea   \label{finalwos}
    \ee
    Since  the  product $\b\tran w_1(N)$ takes place in ${\cal N}_{\!N}$, this means that the 
    neighbors of 
    node $N$ now multiply $w_1(N)$ by their coefficients $\{c_{kN}\}$, for $k \in {\cal N}_N$. That is, we   
    compute the quantity 
    $w_2(k)=c_{kN} w_1(N)$, if  $k \in {\cal N}_N$, otherwise,  we set $w_2(k)=0$. Let the resulting 
    extended vector from this operation   be denoted by $\w_2(1\!:\!N\!-\!1)=\b\tran w_1(N)$. Then 
    (\ref{finalwos}) becomes
    \be 
    \w^o\,=\,  \ba{c}\bar{\L}^{-{\sf T}} \w_2(1\!:\!N\!-\!1) \\ -w_1(N) \ea     \label{finalwos1}
    \ee
    where $\w^o(1\!:\!N\!-\!1)$ is obtained by solving, again, recursively,  
    \be 
    \bar{\L}^{{\sf T}}\w^o(1\!:\!N\!-\!1) = \w_2(1\!:\!N\!-\!1)
    \ee
    
    Table \ref{tab:cgsexdd} summarizes the overall computation of the solution comprising three main passes, 
    viz., $\w_0 
    \rightarrow \w_1 \rightarrow \w_2 \rightarrow \w^o$. Note that $\w_2$ is a sparse vector, and the 
    summations 
    over 
    every row of $\L\tran$ or $\L$ run only within their respective bandwidth.   Observe also that unlike the 
    previous approaches, here the AC is achieved in a sequential manner, with the  
         		first desired mean obtained  in $N$ iterations.  The advantage  in this case is 
    that,   as opposed to the eigenstep algorithm 
    discussed,  
    the back-(forward)\,substitution procedure is well known to be backward stable, and the average consensus 
    is met 
    exactly in $2N$ 
    iterations, regardless of the conditioning of ${\LL}$. 
     
     	\renewcommand{\baselinestretch}{1}
     	\begin{table}[htbp]
     		\begin{center}
     			{\small
     				\begin{tabular}{l} 
     					
     					\noindent {\color{gray!95!white}    \bf

     					\hspace{-.07cm}$_{_{\text{$\blacksquare\!\!\blacksquare\!\!\blacksquare\!\!\blacksquare\!\!\blacksquare\!\!\blacksquare\!\!
     \blacksquare\!\!\blacksquare\!\!\blacksquare\!\!\blacksquare\!\!\blacksquare\!\!\blacksquare\!\!\blacksquare\!\!
     \blacksquare\!\!\blacksquare\!\!\blacksquare\!\!\blacksquare\!\!\blacksquare\!\!\blacksquare\!\!\blacksquare\!\!
     \blacksquare\!\!\blacksquare\!\!\blacksquare\!\!\blacksquare\!\!\blacksquare\!\!\blacksquare\!\!\blacksquare\!\!
     \blacksquare\!\!\blacksquare\!\!\blacksquare\!\!\blacksquare\!\!\blacksquare\!\!\blacksquare\!\!\blacksquare\!\!
     \blacksquare\!\!\blacksquare\!\!\blacksquare\!\!\blacksquare\!\!\blacksquare\!\!\blacksquare\!\!\blacksquare\!\!
     \blacksquare\!\!\blacksquare\!\!\blacksquare\!\!\blacksquare\!\!\!\!\blacksquare\!\!\blacksquare\!\!\blacksquare\!\!\blacksquare\!\!\blacksquare\!\!\blacksquare\!\!\blacksquare\!\!\blacksquare\!\!\blacksquare\!\!\blacksquare\!\!\blacksquare\!\!\blacksquare\!\!\blacksquare\!\!\blacksquare\!\!\blacksquare\!\!\blacksquare\!\!\blacksquare\!\!\blacksquare\!\!\blacksquare\!\!\blacksquare\!\!\blacksquare\!\!\blacksquare$
     					 }}}$}
     					 \\ \\
     %

     					\begin{minipage}{5cm}  
     						\textsl{\phantom{x}for $k=1\!:\!N\!-\!1$,} 
     						\vspace{-.2cm} 
     						\bqn 
     							\;\;\;\;\;\;\;\;\;w_1(k) &\!\!=& \!\!w_0(k)\;- \!\sum_{k>\ell \in 
     						{\cal N}_k} \!\!c_{k\ell}  	w_1(\ell)     
     						\eqn  
     						
     						\vspace{-.5cm}
     						$$
     						 \phantom{x}w_1(N)\;\,=\;\,\gamma \bigg[- w_0(N)\;\,+\sum_{N>\ell \in {\cal 
     						 N}_k} \!\!c_{k\ell} 	 	w_1(\ell)    \biggr]  
     						$$
     						
     						\vspace{-.3cm}
     						\bqn 
     					    \phantom{.}w_2(k) &\!\!=& \left\{ \begin{array}{ll}     \!\!c_{kN} w_1(N), & 
     					    \text{for }k \in 	{\cal N}_N\\ 
     					     \!\!0,  & \text{otherwise} \end{array}   \right.\\
     					 	 \phantom{x}w^o(N)&\!\!=& \!\! -w_1(N)
     					   \eqn
     				 	\textsl{\phantom{x}for $k=N\!-\!1\!:\!1$,} 
     					 	\vspace{-.2cm} 
     					 		$$
     					 		\;\;\;\;\;\;\;\;\;w^o(k) \,=\, w_2(k)\;- \!\sum_{k<\ell \in  	{\cal N}_k} 
     					 		\!\!c_{k\ell}  	w^o(\ell)   
     					 	 $$
     					\end{minipage}	 
     					\phantom{x.} \\ 
     					\vspace{-.1cm}
     					\noindent {\color{gray!95!white}    \bf   
     					\hspace{-.07cm}$_{_{\text{$\blacksquare\!\!\blacksquare\!\!\blacksquare\!\!\blacksquare\!\!\blacksquare\!\!\blacksquare\!\!\blacksquare\!\!\blacksquare\!\!\blacksquare\!\!\blacksquare\!\!\blacksquare\!\!
     \blacksquare\!\!\blacksquare\!\!\blacksquare\!\!\blacksquare\!\!\blacksquare\!\!\blacksquare\!\!\blacksquare\!\!
     \blacksquare\!\!\blacksquare\!\!\blacksquare\!\!\blacksquare\!\!\blacksquare\!\!\blacksquare\!\!\blacksquare\!\!
     \blacksquare\!\!\blacksquare\!\!\blacksquare\!\!\blacksquare\!\!\blacksquare\!\!\blacksquare\!\!\blacksquare\!\!
     \blacksquare\!\!\blacksquare\!\!\blacksquare\!\!\blacksquare\!\!\blacksquare\!\!\blacksquare\!\!\blacksquare\!\!
     \blacksquare\!\!\blacksquare\!\!\blacksquare\!\!\blacksquare\!\!\blacksquare\!\!\blacksquare\!\!\blacksquare\!\!\blacksquare\!\!\blacksquare\!\!\blacksquare\!\!\blacksquare\!\!\blacksquare\!\!\blacksquare\!\!\blacksquare\!\!\blacksquare\!\!\blacksquare\!\!\blacksquare\!\!\blacksquare\!\!\blacksquare\!\!\blacksquare\!\!\blacksquare\!\!\blacksquare\!\!\blacksquare\!\!\blacksquare\!\!\blacksquare\!\!\blacksquare\!\!\blacksquare$
      }}}$}\\\\
     				\end{tabular}
     				\vspace{-.1cm}
     				\caption{\small Exact AC Network back-substitution Algorithm.}\label{tab:cgsexdd}
     			}
     		\end{center}
     		\vspace{-.3cm}
     	\end{table}

     		  Moreover, note that the choice  of $\L$ given by (\ref{LLapl}) is  not unique. There are 
     		alternative ways  one can 
     		think 
     		of in order to minimize, for example, the complexity of the operations in the above recursion. 
     		   For 
     	instance, we can reduce the exchange of information by keeping  a single element   
     	$c_{k\ell}=-1$ below 
     	the main diagonal of $\L$, which still conforms 
     		with the adjacency structure.  As a result, there will be no multiplications involved in the 
     		summations   of Table \ref{tab:cgsexdd},  but only additions. That is, after $2N$ iterations, the 
     		overall 
     		algorithm complexity 
     		amounts to ${\cal O}(2N)$ additions,   as opposed to ${\cal O}(K\!N^2)$ multiplications 
     	that would be necessary for the eigenstep method	to reach the 
     		solution in $K$ iterations via  (\ref{brlsmnl0blms}).  The resulting 
     		recursions in this case  
     		simplify to the ones listed in Table \ref{tab:cgsexdds}, where we further eliminate the 
     		intermediate 
     		variable  $\w_2$, for 
     		compactness.  
     		
     			\renewcommand{\baselinestretch}{1}
     				\begin{table}[htbp]
     					\begin{center}
     						{\small
     							\begin{tabular}{l} 
     								
     								\noindent {\color{gray!95!white}    \bf   
     								\hspace{-.07cm}$_{_{\text{$\blacksquare\!\!\blacksquare\!\!\blacksquare\!\!\blacksquare\!\!\blacksquare\!\!\blacksquare\!\!\blacksquare\!\!\blacksquare\!\!\blacksquare\!\!\blacksquare\!\!\blacksquare\!\!
     								\blacksquare\!\!\blacksquare\!\!\blacksquare\!\!\blacksquare\!\!\blacksquare\!\!\blacksquare\!\!\blacksquare\!\!
     								\blacksquare\!\!\blacksquare\!\!\blacksquare\!\!\blacksquare\!\!\blacksquare\!\!\blacksquare\!\!\blacksquare\!\!
     								\blacksquare\!\!\blacksquare\!\!\blacksquare\!\!\blacksquare\!\!\blacksquare\!\!\blacksquare\!\!\blacksquare\!\!
     								\blacksquare\!\!\blacksquare\!\!\blacksquare\!\!\blacksquare\!\!\blacksquare\!\!\blacksquare\!\!\blacksquare\!\!\blacksquare\!\!\blacksquare\!\!\blacksquare\!\!
     									\blacksquare\!\!	\blacksquare\!\!	
     									\blacksquare\!\!\blacksquare\!\!	
     									\blacksquare\!\!	
     									\blacksquare\!\!\blacksquare\!\!\blacksquare\!\!\blacksquare\!\!\blacksquare\!\!\blacksquare\!\!\blacksquare\!\!\blacksquare\!\!\blacksquare\!\!\blacksquare\!\!\blacksquare\!\!\blacksquare\!\!\blacksquare\!\!\blacksquare\!\!\blacksquare\!\!\blacksquare\!\!\blacksquare\!\!\blacksquare\!\!\blacksquare\!\!\blacksquare\!\!\blacksquare\!\!\blacksquare\!\!\blacksquare\!\!\blacksquare\!\!\blacksquare\!\!\blacksquare\!\!\blacksquare\!\!\blacksquare\!\!\blacksquare\!\!\blacksquare\!\!\blacksquare\!\!\blacksquare\!\!\blacksquare\!\!\blacksquare$
     								 }}}$}\\\\
     			%

     								\begin{minipage}{5cm}  
     							 	\vspace{-.5cm} 
     									\bqn 
     									 \;\;w_1(k) &\!\!=& \!\!w_0(k)\;+ \!\sum_{k>\ell \in 
     									{\cal N}_k} \!\! 	w_1(\ell)    \;, \;\;\;k=1,2,\ldots,N\!-\!1 
     									\eqn  
     									\vspace{-.1cm}
     									$$
     									 \phantom{x}w^o(N)\;\,=\;\, 
     									 w_0(N)\;\,+\sum_{N>\ell \in 
     									 {\cal 
     									 N}_k} \!\! 	 	w_1(\ell)     
     									$$
     									\vspace{-.3cm}
     									\bqn 
     								    \phantom{x}w_1(k) &\!\!=& \left\{ \begin{array}{ll}     \!\!   
     								    \gamma  w^o(N), 
     								    & 
     								    \text{for }k \in 	{\cal N}_N\\ 
     								     \!\!0,  & \text{otherwise} \end{array}   \right.
     								   \eqn
     								 	\vspace{-.2cm} 
     								 		$$
     								 	 \;\;w^o(k) \,=\, w_1(k)\;+ \;	w^o(\ell)   \;,\;\;k<\ell \in  	
     								 		{\cal N}_k\;,\;   k=N\!-\!1,\ldots,2,\!1
     								 	 $$
     								\end{minipage}	 
     								\phantom{x.} \\ 
     								\vspace{-.1cm}
     								\noindent {\color{gray!95!white}    \bf   
     								\hspace{-.07cm}$_{_{\text{$\blacksquare\!\!\blacksquare\!\!\blacksquare\!\!\blacksquare\!\!\blacksquare\!\!\blacksquare\!\!\blacksquare\!\!\blacksquare\!\!\blacksquare\!\!\blacksquare\!\!\blacksquare\!\!
     							     	\blacksquare\!\!\blacksquare\!\!\blacksquare\!\!\blacksquare\!\!\blacksquare\!\!\blacksquare\!\!\blacksquare\!\!
     							    	\blacksquare\!\!\blacksquare\!\!\blacksquare\!\!\blacksquare\!\!\blacksquare\!\!\blacksquare\!\!\blacksquare\!\!
     							     	\blacksquare\!\!\blacksquare\!\!\blacksquare\!\!\blacksquare\!\!\blacksquare\!\!\blacksquare\!\!\blacksquare\!\!
     							     	\blacksquare\!\!\blacksquare\!\!\blacksquare\!\!\blacksquare\!\!\blacksquare\!\!\blacksquare\!\!\blacksquare\!\!\blacksquare\!\!\blacksquare\!\!\blacksquare\!\!
     							     	\blacksquare\!\!		\blacksquare\!\!	
     							     	\blacksquare\!\!\blacksquare\!\!	\blacksquare\!\!	
     							     	\blacksquare\!\!\blacksquare\!\!\blacksquare\!\!\blacksquare\!\!\blacksquare\!\!\blacksquare\!\!\blacksquare\!\!\blacksquare\!\!\blacksquare\!\!\blacksquare\!\!\blacksquare\!\!\blacksquare\!\!\blacksquare\!\!\blacksquare\!\!\blacksquare\!\!\blacksquare\!\!\blacksquare\!\!\blacksquare\!\!\blacksquare\!\!\blacksquare\!\!\blacksquare\!\!\blacksquare\!\!\blacksquare\!\!\blacksquare\!\!\blacksquare\!\!\blacksquare\!\!\blacksquare\!\!\blacksquare\!\!\blacksquare\!\!\blacksquare\!\!\blacksquare\!\!\blacksquare\!\!\blacksquare\!\!\blacksquare$
     							     								 }}}$}\\\\
     							\end{tabular}
     							\vspace{-.3cm}
     							\caption{\small Low complexity, Exact AC algorithm of Table 
     							I.}\label{tab:cgsexdds}
     						}
     					\end{center}
     					\vspace{-.2cm}
     				\end{table}

     	\noindent  {\sfb 2)  \underline{Graph Filter implementation}}. \vspace{.1cm}
     	
     		Consider the product $ 
     	\w_1(1\!\!:\!\!N\!-\!1)= \bar{\L}^{-1}   
     		\w_0(1\!\!:\!\!N\!-\!1) $ 
     		in 
     		(\ref{expswow0}). As opposed to solving 
     		it via back-substitution, we can implement the   $(N\!-\!1) \times(N\!-\!1) $ inverse 
     		$\bar{\L}^{-1}$ in a 
     		finite number of iterations as follows. Let $\bar{\L}=(\I+\L_o)$, where $\L_o$ is   strictly 
     		lower 
     		triangular. Since  $\L_o$ is a nilpotent matrix, the following expansion into a  finite number of 
     		terms   
     		can be 
     		easily 
     		verified: 
     		\be 
     		(\I+\L_o)^{-1}\,=\,\I \,+\,\sum_{l=1}^{N-2} (-1)^l \L_o^l \;\define\; f(\L_o) \label{polymat}
     		\ee
     	 where the polynomial $f(\L_o)$ can be interpreted as a matrix FIR filter in terms of the 
     		``shift'' operator  $\L_o$. The former has been generally referred to in graph theory as  a   
     		{\em 
     		graph 
     		filter}\,\cite{Moura2013}.  
     		  
     		  Now, let  $\z_{N-2} \define \w_1(1\!:\!N\!-\!1)$ and 
     		  $\z_0=\w_0(1\!\!:\!\!N\!-\!1)$, for 
     		  compactness of notation. Then, 
     		  \be
     		   \z_{N-2}= \z_0 - \L_o\z_0 + \L_o^2\z_0 - 
     		   \L_o^3\z_0 +\;\ldots  \;  (-1)^{N-2}\L_o^{N-2} \label{reczl}
     		   \ee
     		   which  naturally yields the recursion, for $l=1,2,\ldots,N-2$,
     		   \be 
     		   \z_l\,=\,\z_0 -\L_o\z_{l-1}  \;,\;\;\;\;\;\z_0=\w_0(1\!:\!N\!-\!1) \label{recyl}
     		   \ee 
     		   Similarly, the product $\y_{N-2} \define \w^o(1\!:\!N\!-\!1)=\bar{\L}^{-{\sf T}}\b\tran 
     		   \w_1(N)$  in 
     		   (\ref{finalwos}),   can be computed   analogously  as 
     		   \be 
     		   	   \y_l\,=\, \y_0 - \L_o\tran \y_{l-1}  \;,\;\;\;\;\;\;\;\;\; \y_0=\b\tran \w_1(N)
     		   	\ee
     		  
     		  Note that the 
     	 structure of $\L_o$ may be such that $\L_o^p$  vanishes for $p \!<\!\!< \!N-2$. In addition,  the 
     	  exchange of information  can again be   minimized 
     	 by reducing the number of elements in every column of $\L$ down to two, i.e., to $\{1,-1\}$. 
     	 We   illustrate this fact  using the  scenario   of the matrix in Fig.\;\ref{fig.Adjmat1}, in which 
     	 case 
     	 we 
     	 find that $\L_o^{23}=0$. As a result, the  recursions  (\ref{reczl}) and  (\ref{recyl}) along with 
     	 the 
     	 effective matrix multiplication with $\b\tran \b$, amount  to $47$ iterations {\em matrixwise}. This 
     	 can 
     	 again be 
     	 compared  to    
     	 the 
     	 number of 
     	 iterations   required by the eigenstep recursion. That is, by 
     	 computing the eigenvalues 
     	 of $\L \L\tran$, we verified that the eigenvalue at 1 has multiplicity 6, yielding $N-6=47$, which 
     	 is 
     	 the same 
     	 required 
     	 by the eigenstep  algorithm. This is  expected, since  the number of iterations governed   by its 
     	 eigenstructure in the latter should be somehow translated into the number multiplying matrices 
     	 appearing in the 
     	 implementation of (\ref{polymat}).  	 Indeed, observe that the overall computation in this example 
     	 is a product of $N-6$ matrices of the 
     	 form 
     	 (\ref{prodsetW}), since $f(x)$ in (\ref{polymat}) admits a standard factorization of the form 
     		\be 
     		f(\L_o)\,=\,(-1)^{N-2}\prod_{k=1}^{N-2}(\L_o-\alpha_k\I) \label{polymatfact}
     		\ee
     	where $\alpha_k$ are the  roots of the polynomial $f(x)$\,(which appear as 
     	complex conjugates on the unit circle).   That is, by defining the matrices $\L_{\alpha_k} \define 
     	\L_o-\alpha_k\I$,  
     	(\ref{prodmatswo}) can be written as
     	\be 
     	\w^o\,=\,\gamma \left(\prod_{k=1}^{23} \ba{cc}  \!\!\L_{\alpha_k}\tran & \0\\ \0 & 1\ea\right) 
     	\b_1\tran 
     	\b_1\left( 
     	\prod_{k=1}^{23} \ba{cc}  
     	 \!\!\L_{\alpha_k}  & \0\\ \0 & 1\ea\right) \w_0 \label{2ndmeth}
     	\ee
     	Since $\L_0$ with minimized number of elements has $N\!-\!2$ components, each  matrix 
     	iteration  
     	involves only ${\cal O}(N)$ 
     	additions, so that the overall computation will require ${\cal O}(N^2)$ 
     	additions.\!\footnote{\noindent 
     The  formula (\ref{2ndmeth}) could be extended to directed graphs by 
     	starting   from 
     	   \be 
     	   \w^o\,=\,\Big[\I -    \C_2\tran  \big( \C_1   \C_2\tran \big)^{\!-1}  \C_1\Big]\w_0
     	    \ee
     	    and   selecting $\{\C_1,\C_2\}$ as  upper triangular matrices of the form 
     	    (\ref{Cpartix}) as
     	      \be 
     	       \C_1\,=\,\L_1\tran\,=\, \Big[ \xoverline{\L}_1\tran \;\;\; \b\tran\Big] \,,\;\;\;\; 
     	       \C_2\,=\,\L_2\tran\,=\,\Big[ 
     	       \xoverline{\L}_2\tran \;\;\; \g\tran\Big] \label{CpartixL}
     	       \ee
     	       The construction of their corresponding Laplacian-like structures, however, is not 
     	       straightforward, since a block tridiagonal structure should result from some 
     	       permutation strategy in first place. If this is succeeded, the formulas would follow   
     	       similarly to
     	       the ones of the symmetric case, where now
     	       \be
     	      \w^o 
     	     \,=\,\gamma \ba{cc} \! \xoverline{\L}_1^{-{\sf T}} & \0 \\ \0 & 1 
     	               \ea 
     	              \! \ba{c}   \!\b\tran \\   -1\ea    \! \ba{cc} \!\g & -1 
     	              \ea  
     	              \! \ba{cc} \!\! 
     	             \xoverline{\L}_2^{-1}  & \0 \\ \0 &  1 \ea   \w_0 \nonumber 
     	         \ee}

\section{AC in Diffusion Adaptive Schemes} \label{ACdiffu}

 The fact that the exact solution (\ref{brlsmnl0002}) intrinsically 
 performs averaging, suggests that some related problems of   diffusion adaptation over networks should make 
 use of fast AC,  or, for the same matter, a reduced complexity algorithm in their designs.  That is,  assume 
 that  each node $k$  
 of 
 a network of 
 agents receives streaming 
  data $\{d_k(n),\u_{k,n}\}$,   related via a simple   linear regression model of the form
  $d_{k}(n) = \u_{k,n}\w_{k}^o + v_{k}(n)$,
  where $n$ is the time index, $\w_k^o$ is an unknown  local parameter of size $M\times 1$, $\u_{k,n}$ is a 
  regression (row)\,vector of size $1\times M$, and        $v_k(n)$ is an  additive,  zero-mean white 
  noise.   
  The latter is  assumed to be time and spatially uncorrelated with other data.  Then, a global linear 
  model for 
  the 
  agents is 
  commonly 
described by means of  the following extended definitions\,:
  \bq
  \WW^o  &=& {\sf col}\,\{\w^o_1 ,\w^o_2  ,\ldots,\w^o_N\} \label{eqentriesop2}\\
  \d_n &=& {\sf col}\,\{d_{1}(n),   d_{2}(n),   \ldots , d_{N}(n)\} \label{Ymatrix00} \\
   \U_{n}&=&\bdiag\left\{ \u_{1,n} ,\u_{2,n} ,\ldots,\u_{N,n} \right\}  \label{Ymatrix} \\
   \v_n  &=&{\sf col}\,\{ v_{1}(n) , v_{2}(n),\ldots, v_{N}(n)\} \label{netnoisek}
  \eq
  so that    
  \be
  \d_n  \;=\; \U_{\!n}\WW^o + \v_n \label{linmodelwp2}
  \ee
Let   $\otimes$  denote   the {\it Kronecker} product,  and $\Ex$,   the expectation operator.  
By defining $\Ccal=\C \otimes  \I_{\!M}$, one can seek the solution, say, ${\WW}_{n}$, to the  following  
network constrained  LMS problem: 
      \be
      \min\limits_{\WW} \; \Ex \!\left\| \d_{n} -  \U_{n}\WW \right\|^2 \;\;\,\,\,\,\,\,\;\text{s.t. }\;\,  
      \Ccal \WW=\0  \label{LCreggenc000b}
      \ee 
   After completing squares, this cost becomes equivalent to  (\ref{compfcccwo}), where we can
     select, as the prior guess for the solution, 
     an     unconstrained  LMS   update  from the previous estimate, ${\WW}_{n-1}$. That is, given some 
     step-size 
               $\mu$, and by setting $\w_0 
     \longleftarrow 
     {\WW}^{u}_n$,   with $\Pb_0\longleftarrow \I$, we compute: 
      \bq
     {\WW}^{u}_{n}&=&{\WW}_{n-1} +  \mu \U_{n}^\star  (\d_n -  \U_n\WW_{n-1})\label{brlccccnlbx} \\
      {\WW}_{n} &=&\big[\I\!-\!\Ccal^\star(\Ccal\Ccal^ \star)^{-1} \Ccal \big] {\WW}^{u}_n   
      \label{brlccccnlb}
      \eq

Expression (\ref{brlccccnlb})   suggests that if one is to minimize the complexity of this step, this has to 
take place as  
 efficiently  as possible.   In other words,    if averaging 
 is to 
 be accomplished by several consecutive  matrix 
 iterations 
 between every self-learning update (\ref{brlccccnlbx}),  the combination matrix $\Ccal$ should be designed 
 to 
 achieve fast averaging,   
 and 
 not 
 necessarily as any arbitrary (and more commonly) stochastic matrix. This can be contrasted, e.g., with the 
 CG 
 method of 
 \cite{mer22B},   
 which in the 
 distributed scenario will not be exact, it will be more complex, and still dependent on the eigenvalue 
 spread of 
 the 
 combination matrix.
 \vspace{.1cm} 
     	
     	 \section{Conclusions}
    
    We  have revisited the AC strategy in a network of distributed agents as a linearly-constrained  
    problem. 
    We showed that when the solution is implemented via  eigensteps,    numerical 
    errors 
    can accumulate  to the point where 
    the consensus cannot be met.  
    The proposed alternatives  simply borrow from   well known standalone 
    counterparts by capitalizing on the fact that, in a descentralized mode, both combination matrix 
     and its
     eigenpairs are fully known.  The recursions in this case   do not involve    online computation of 
     global 
     parameters, and offer one 
    possible 
    stable solution to the problem, depending on its size. This includes the standard  power iterations and 
    the 
    eigenstep method, both  extended to a  broader case of directed networks. 
    
   From a different perspective, we showed that the exact AC can be reached in 
   the same number of iterations that would be required by the eigensteps   via (\ref{2ndmeth}) with 
   ${\cal O}(N^2)$ additions, while complexity can be 
    reduced  to $2N$ additions with the algorithm of Table\;II. The key aspect in the latter is the fact 
    that,  by properly 
    renumbering   nodes, we can draw a path of accumulations from node 1  to  node $N$,  
  yielding the desired means sequentially, as we trace our way back to node 1. Unlike the previous iterative 
  methods, this is independent of the conditioning of the underlying constraint matrix.
  
     We hasten to add that, although 
  iterative 
  methods 
  are 
  not exact,   
  they can potentially  yield estimates that are 
      closer 
      to  the solution compared to an exact algorithm, depending on how the time-scale is measured. This can 
      be 
      more relevant to adaptat-and-diffuse schemes  like (\ref{brlccccnlbx})-(\ref{brlccccnlb}),  
      when a 
      limited 
      number of 
      iterations is allowed for approximating (\ref{brlccccnlb}).

\appendices
\section{} \label{ApA}
  
  Table\;\ref{tab:combm} lists  the nodes/edges information for the directed graph 
  corresponding to
  Fig.\,\ref{fig.net_asymm_20}\,($\A_o$), in a standard form. The middle column refers to the 
  indexes of   nonzero row elements of $\A_o$, with  the respective edges listed on the right side. 
    	\begin{table}[htbp]
       \begin{center}
    $$
     \tiny {\begin{array}{ccr} \text{\scriptsize {\bf row}} & \text{\scriptsize {\bf  column}} & 
     \text{\scriptsize  {\bf edge}} 
   \\  \hhline{===}   1 & 1 & 
    0.3187\\ 
    1 & 4 & 
    -0.0243\\ 
    1 & 14 & 0.2208\\ 1 & 15 & 0.0766\\ 2 & 2 & 
    0.1773\\ 2 & 14 & 0.0487\\ 2 & 15 & 0.1334\\ 2 & 17 & 0.2641\\ 3 & 3 & 0.00041106\\ 3 & 4 & 0.1153\\ 3 & 
    5 
    & -0.0289\\ 3 & 13 & 0.0869\\ 3 & 15 & -0.0115\\ 3 & 19 & 0.1878\\ 4 & 1 & 0.1410\\ 4 & 3 & 0.1216\\ 4 & 
    4 
    & 0.3227\\ 4 & 5 & 0.1237\\ 4 & 15 & 0.1266\\ 5 & 3 & 0.3486\\ 5 & 4 & 0.0323\\ 5 & 5 & 0.0134\\ 5 & 10 & 
    0.3539\\ 5 & 13 & 0.0554\\ 5 & 14 & -0.0049\\ 5 & 15 & 0.1562\\ 6 & 2 & 0.0762\\ 6 & 6 & 0.0323\\ 6 & 8 & 
    0.0427\\ 6 & 9 & 0.1354\\ 7 & 7 & 0.1218\\ 7 & 8 & 0.3948\\ 7 & 10 & 0.2981\\ 7 & 14 & -0.0034\\ 7 & 16 & 
    0.1269\\ 8 & 2 & 0.2228\\ 8 & 6 & 0.2379\\ 8 & 8 & 0.2200\\ 8 & 9 & -0.0143\\ 8 & 10 & 0.0501\\ 9 & 2 & 
    0.2137\\ 9 & 6 & 0.0765\\    9 & 8 & 0.0247\\ 9 & 9 & 0.2090\\   10 & 5 & -0.0182\\  10 & 8 & 0.0259\\  
    \hhline{---} 
        \end{array} \vline   \begin{array}{ccr} 
                \text{\scriptsize {\bf row}} & \text{\scriptsize {\bf column}} & 
                     \text{\scriptsize  {\bf edge}} 
                   \\ 
                    \hhline{===} 10 & 10 & 
    0.1709\\ 10 & 13 & 0.1681\\ 10 & 16 & 0.0463\\ 10 & 17 & 0.2847\\ 11 & 11 & 0.2769\\ 11 & 17 & 0.2550\\ 
    12 
    & 7 & 0.0766\\ 12 & 12 & 0.1096\\ 12 & 14 & 0.2778\\ 12 & 18 & 0.2140\\ 13 & 3 & 0.0935\\ 13 & 5 & 
    0.1750\\ 
    13 & 10 & -0.0015\\ 13 & 13 & 0.1726\\ 14 & 1 & 0.2030\\ 14 & 5 & -0.0152\\ 14 & 7 & 0.1430\\ 14 & 12 & 
    0.3554\\ 14 & 14 & 0.2865\\ 14 & 15 & 0.0802\\ 14 & 18 & 0.3390\\ 14 & 20 & -0.0379\\ 15 & 1 & -0.0246\\ 
    15 
    & 10 & 0.1201\\ 15 & 11 & 0.3375\\ 15 & 12 & 0.3285\\ 15 & 14 & 0.3454\\ 15 & 15 & 0.3152\\ 16 & 10 & 
    0.1246\\ 16 & 16 & 0.1955\\ 16 & 17 & 0.1620\\ 17 & 5 & -0.0178\\ 17 & 6 & 0.4003\\ 17 & 16 & 0.1499\\ 17 
    & 
    20 & 0.2325\\ 18 & 7 & 0.3575\\ 18 & 12 & 0.2434\\ 18 & 14 & 0.0322\\ 18 & 18 & 0.0199\\ 18 & 20 & 
    0.0961\\ 
    19 & 3 & 0.2130\\ 19 & 19 & 0.0263\\ 20 & 7 & 0.3310\\ 20 & 14 & 0.2856\\ 20 & 18 & 0.0366\\ 20 & 20 & 
    0.1546 \\  \hhline{---} \end{array}
    }
   $$
   \caption{\small Directed Graph information for the network of 
   Fig.\,\ref{fig.net_asymm_20}.}\label{tab:combm}
    \end{center}
    \end{table}

    \section*{Acknowledgement}
   
   The author would like to thank Prof.\;Vitor\;H.\;Nascimento, from the 
   University 
   of  
   São Paulo\,(USP), Brazil,   for his 
   valuable comments   during the   preparation of this manuscript.

{\small

}

\begin{IEEEbiographynophoto}{Ricardo Merched}
received the BS  and  MS degrees in 1995 and 1997  from  the Federal 
University of Rio de Janeiro\,(UFRJ), Brazil, and the Ph.D. degree   from the University 
of California, Los Angeles (UCLA), 
in 
2001, all in Electrical Engineering.  He became Professor at the Departm. of Electronic  
Engineering, UFRJ, in 2002. He was 
a 
Visiting Scholar with UCI, the University of California, Irvine\,(2006-2007),  the  University Graduate 
Center 
in Oslo\,(UniK)\,(2010), and the  École Polytechnique Fédérale de Lausanne\,(EPFL) in 2018-2019. 
He was Associate Editor of the {\em IEEE 
  TRANSACTIONS ON CIRCUITS AND SYSTEMS I}, the {\em IEEE TRANSACTIONS ON SIGNAL PROCESSING LETTERS}, and 
 the 
EURASIP,  the
{\em  European Journal on Advances in Signal Processing}. His current  interests include machine 
learning,   efficient DSP techniques for MIMO   communications, and  distributed processing over 
complex networks.
\end{IEEEbiographynophoto}


\end{document}